\documentclass[reqno]{amsart}     
\usepackage{enumerate,amsmath}     
\usepackage{pb-diagram}  
\usepackage{pstricks,pst-node,pst-coil,pst-plot}   
     
\theoremstyle{plain}      
\newtheorem{step}{Step} 
\newtheorem{thm}{Theorem}[section]     
\newtheorem{theorem}[thm]{Theorem}     
\newtheorem{cor}[thm]{Corollary}     
     
\newtheorem{lem}[thm]{Lemma}     
\newtheorem{lemma}[thm]{Lemma}     
\newtheorem{prop}[thm]{Proposition}     
\newtheorem{claim}[thm]{Claim}     
\theoremstyle{remark}      
\newtheorem{rem}[thm]{Remark}     
\theoremstyle{definition}      
     
\newtheorem{definition}[thm]{Definition}     

\def\al{{\alpha}}         
         
\def\de{{\delta}}         
         
\def\om{{\omega}}         
\def\Om{{\Omega}}         
\def\la{{\lambda}}

\def\ep{{\varepsilon}}

\def\phi{{\varphi}}

\let\pa\partial     
\let\na\nabla     
     
\DeclareMathAlphabet{\doba}{U}{msb}{m}{n}

\gdef\mN{\doba{N}}

\gdef\mR{\doba{R}}         
\gdef\mS{\doba{S}}         
         
\gdef\mZ{\doba{Z}}

\def\w{\overline{w}}
\def\v{\overline{v}}
   
\def\Vol{{\mathop{\rm Vol}}}     
\def\Scal{{\mathop{\rm Scal}}}     
     
\def\Spec{{\mathop{\rm Spec}}}     
\let\ti\tilde   
\let\ol\overline
\def\eref#1{{\rm (\ref{#1})}}   

\let\<\langle 
\let\>\rangle 
\def\spann{\mathop{\rm span}}
\def\Vzero{V\setminus\{0\}}
\let\weakto\rightharpoonup

\def\theobull{\ \ $\bullet$}
\let\ss\scriptstyle

\def\Gr#1#2#3{{\rm Gr}_{#1}^{#2}(#3)}

\begin{document}     
\title{The second Yamabe invariant}

\maketitle     
\begin{center}
\sc B. Ammann and E. Humbert\footnote{bernd.ammann@gmx.net, humbert@iecn.u-nancy.fr}
\end{center}     

\begin{abstract}     
Let $(M,g)$ be a compact Riemannian manifold of dimension $n \geq 3$. We
define the second Yamabe invariant  as the infimum of the second eigenvalue
of the Yamabe operator over the metrics
conformal to~$g$ and of volume~$1$. We study when it is
attained. As an application, we find nodal solutions of the Yamabe equation.
\end{abstract}     
  
\begin{center}
\today
\end{center}   


\tableofcontents

     
\section{Introduction} 
Let $(M,g)$ be an $n$-dimensional compact Riemannian manifold ($n \geq 3$). 
In \cite{yamabe:60} Yamabe attempted to show that there is a metric           
$\tilde{g}$  conformal to $g$ such that the scalar curvature            
$S_{\tilde{g}}$ of ${\tilde{g}}$ is            
constant. However, Trudinger \cite{trudinger:68}
realized that Yamabe's            
proof contained a serious gap.           
The problem is now solved, but it took a very long time to find the good
approach. The problem of finding a metric ${\tilde{g}}$
with constant scalar curvature in the conformal class $[g]$          
is called the Yamabe problem. 
The first step towards a rigorous solution of this problem was achieved by 
Trudinger~\cite{trudinger:68} who was able to repair the
gap of Yamabe's article in the case that the scalar curvature of $g$ is 
non-positive.
Eight years later,   
Aubin \cite{aubin:76} solved the problem  for arbitrary non locally
conformally flat manifolds of dimension  $n \geq 6$. 
The problem was completely solved 
another eight years later in an article of Schoen \cite{schoen:84} in which 
the proof was reduced to the positive-mass theorem which had previously been 
proved by Schoen and Yau \cite{schoen.yau:79a,schoen.yau:88}. The reader can refer to
\cite{lee.parker:87}, \cite{aubin:76} or \cite{hebey:97} 
for more information on this
subject. The method to solve the Yamabe problem was the following. Let $u \in
C^{\infty}(M)$, $u>0$ be a smooth function and $\tilde{g} = u^{N-2} g$
where $N= \frac{2n}{n-2}$. Then, multiplying $u$ by a constant,
the following equation is satisfied:
\begin{eqnarray*} 
L_g (u)  =  S_{\tilde{g}} |u|^{N-2} u. 
\end{eqnarray*}
where 
$$L_g= c_n \Delta_g +S_g =   {4(n-1) \over n-2} \Delta_g + S_g$$
is called the Yamabe operator. As a consequence, solving the 
Yamabe problem is equivalent
to finding a positive smooth solution $u$ of   
\begin{eqnarray} \label{eqyam} 
L_g (u)  = C_0   |u|^{N-2} u. 
\end{eqnarray}
where $C_0$ is a constant. In order to obtain solutions of this equation
Yamabe defined the 
quantity
$$\mu(M,g)= \inf_{u \not= 0, u \in C^{\infty}(M)} Y(u)$$ 
where 
$$Y(u)= \frac{\int_M c_n |\nabla u|^2 + S_g u^2\,dv_g}{{ \left( 
\int_M |u|^N \,dv_g \right)}^{\frac{2}{N}}}.$$
Nowadays, $\mu(M,g)$  is called the 
\emph{Yamabe invariant}, and $Y$ the   
\emph{Yamabe functional}. Writing the Euler-Lagrange equation associated to
$Y$, we see that there exists a one to one correspondence between critical
points of $Y$ and solutions of equation (\ref{eqyam}). In particular, if $u$
is a positive smooth function
such that $Y(u)= \mu(M,g)$, then $u$ is a solution of (\ref{eqyam}) and 
$\tilde{g}= u^{N-2} g$ is the desired metric of constant scalar curvature. 
The key point of
the resolution of the Yamabe problem is the following theorem due to 
Aubin \cite{aubin:76}. In the theorem and in the whole article, $\mS^N$
will always denote the sphere $S^n$ with the standard Riemannian structure.
\\

\begin{theorem} \label{aubin}
Let $(M,g)$ be a compact Riemannian manifold of dimension $n \geq 3$. If 
$\mu(M,g) < \mu(\mS^n)$, then there exists a positive smooth
function $u$ such that $Y(u)= \mu(M,g)$.
\end{theorem}

This strict inequality
is used to show that a minimizing sequence does not concentrate
in any point. Aubin
\cite{aubin:76} and
Schoen \cite{schoen:84} proved the following.

\begin{theorem} 
Let $(M,g)$ be a compact Riemannian manifold of dimension $n \geq 3$. Then 
$\mu(M,g) \leq \mu(\mS^n)=  n(n-1) \om_n^{\frac{2}{n}}$ 
where $\om_n$ stands for
the volume of the standard sphere $S^n$. Moreover, we have equality in this
inequality if and only if $(M,g)$ is conformally
diffeomorphic to the sphere.
\end{theorem}
These theorems solves the Yamabe problem. \\

In this paper, we introduce and study an invariant 
that we will call the \emph{second Yamabe invariant}. 
It is well known that the operator $L_g$ has discrete spectrum 
$$\Spec(L_g)= \{ \la_1(g), \la_2(g), \cdots \}$$
where the eigenvalues
$$\la_1(g) < \la_2(g) \leq \la_3(g) \leq \cdots \leq \la_k(g) \cdots \to +\infty$$
appear with their multiplicities. 
The variational characterization of $\la_1(g)$ is given by 
$$\la_1(g)= \inf_{u \not= 0, u \in C^{\infty}(M)} 
\frac{\int_M c_n |\nabla u|^2 + S_g u^2\,dv_g}{\int_M |u|^2 \,dv_g}.$$
Let $[g]$ be the conformal class of $g$.  
Assume now that the Yamabe invariant $\mu (M,g) \geq 0$. 
It is easy to check that  
$$\mu(M,g) = \inf_{\tilde{g} \in [g]} \la_1(\tilde{g})
  \Vol(M,\tilde{g})^{\frac{2}{n}},$$
where $[g]$ is the conformal class of $g$.
We then  enlarge this definition.

\begin{definition} Let $k \in \mN^*$. Then, the $k^{th}$ Yamabe invariant is
  defined by 
$$\mu_k(M,g) = \inf_{\tilde{g} \in [g]} \la_k(\tilde{g})
  \Vol(M,\tilde{g})^{\frac{2}{n}}.$$
\end{definition}
 
With these notations, $\mu_1(M,g)$ equals to Yamabe invariant $\mu(M,g)$ 
in the case $\mu(M,g) \geq 0$, and $\mu_1(M,g)=-\infty$ in the case 
$\mu(M,g) < 0$. 

The goal of this article is to study the second Yamabe invariant
$\mu_2(M,g)$ for manifolds whose Yamabe invariant in the case  
$\mu(M,g) \geq 0$.
As explained in Section \ref{negative}, the most interesting
case is when $\mu(M,g)>0$.
In particular, we discuss whether 
$\mu_2(M,g)$ is attained. This question is
discussed in Subsection~\ref{attaine}. 
In particular,  Proposition~\ref{notattained} asserts that
contrary to the standard Yamabe invariant, $\mu_2(M,g)$ cannot  be attained
by a metric if $M$ is connected. In other words, 
there does not exist
$\tilde{g} \in [g]$ such that $\mu_2(M,g) = \la_2(\tilde{g})
\Vol(M,\tilde{g})^{\frac{2}{n}}$. 
In order to find minimizers, we enlarge the conformal class $[g]$
to what we call the class of \emph{generalized metrics} conformal to $g$.
A generalized metric is a ``metric'' of the form 
$\tilde{g} = u^{N-2} g$, where $u$ is no longer necessarily 
positive and smooth, but 
$u \in L^{N}(M)$, $u\geq 0$, $u\not\equiv 0$.
The definitions of $\la_2(\ti g)$ and of $\Vol(M,\ti g)$ 
can be extended to generalized
metrics (see section 3). Then, we are able to prove
the following result:

\begin{theorem} \label{attain}
Let $(M,g)$ be a compact Riemannian manifold of dimension $n \geq 3$ whose
Yamabe invariant is non-negative. Then,
$\mu_2(M,g)$ is attained by a generalized metric in the following cases:
\begin{enumerate}[\theobull]
\item  $\mu_1(M,g) > 0 $ and $\mu_2(M,g) < \left[\mu_1(M,g)^{\frac{n}{2}} +
\mu_1(\mS^n))^{\frac{n}{2}} \right]^{\frac{2}{n}} $; 
\medskip
\item  $\mu_1(M,g) = 0 $ and $\mu_2(M,g) < \mu_1(\mS^n)$

\end{enumerate}
where  $\mu_1(\mS^n)= n(n-1) \om_n^{\frac{2}{n}} $ is the Yamabe invariant of the standard sphere.
\end{theorem}

The result we obtain in the case $ \mu_1(M,g) = 0 $ is not
surprising. Indeed, when $\mu_2(M,g) < \mu_1(\mS^n)$, 
Aubin's methods \cite{aubin:76} can be adapted here
and 
allow to avoid concentration of
minimizing sequences.
However, when  $\mu_1(M,g) > 0 $ and $\mu_2(M,g) < \left[\mu_1(M,g)^{\frac{n}{2}} +
\mu_1(\mS^n))^{\frac{n}{2}} \right]^{\frac{2}{n}} $, the result is much
more difficult to obtain (see Subsection~\ref{sectatt}). 
A second result is to find explicit examples for which 
the assumptions of Theorem \ref{attain} are satisfied. The method consists
in finding an appropriate couple of test functions. 
\begin{theorem} \label{condi}
The assumptions of Theorem  \ref{attain} are satisfied in the following cases:
\begin{enumerate}[\theobull]
\item $\mu_1(M,g)>0$, $(M,g)$ is not locally conformally flat and $n \geq
11$;
\medskip
\item $\mu_1(M,g)=0$, $(M,g)$ is not locally conformally flat and $n \geq
9$.
\end{enumerate}

\end{theorem}
One of our motivations is to find solutions of the Yamabe
equation (\ref{eqyam}) with alternating sign, i.e. positive and negative
values. 
If $M$ is connected, alternating sign implies
that the zero set $u^{-1}(0)$ of $u$ is not empty. In the following we 
will use the standard definition to call the zero set  $u^{-1}(0)$ of a function 
$u$ the \emph{nodal set of $u$}. A solution with a non-empty nodal set
is usually called a \emph{nodal solution}. If $M$ is connected, then the 
maximum principle implies that a solution of the Yamabe equation is nodal 
if and only if it has alternating sign.
They are called \emph{nodal solutions} 
of the Yamabe equation. The articles \cite{hebeyvaugon:94}, 
\cite{djadli.jourdain:02},
\cite{jourdain:99}, \cite{holcman:99} prove existence of nodal solutions
under symmetry assumptions or under some assumptions which allow to use
Aubin's methods, as in Theorem~\ref{attain} when  $\mu_1(M,g) = 0 $ and
$\mu_2(M,g) < \mu_1(\mS^n)$. If $\mu(M,g) \leq 0$, another method is given
in Section \ref{negative}. 
The method we use here is completely
different and we obtain solutions on a large class of manifolds. 
In particular, to our knowledge, there
is no work which leads to the existence of such solutions if the Yamabe
invariant is positive and if $(M,g)$ is not conformally equivalent
to the round sphere. The result we obtain  is the
following:

\begin{theorem} \label{eulerequ}
 Let $(M,g)$ be a compact Riemannian manifold of dimension $n \geq 3$.
Assume that $\mu_2(M,g)$ is attained by a generalized metric $u^{N-2}g$
where $u \in L^N(M)$, $u \geq 0$ and $u \not\equiv  0$. Let $\Om $ be the nodal set of $u$. Then, there
exists a nodal solution $w \in C^{\infty}(M \setminus \Om) \cap C^{3,\alpha}(M)$ (
$\alpha \leq N-2$) of equation (\ref{eqyam}) such that $|w| =u $. 
\end{theorem}

 A corollary of Theorems \ref{attain}, \ref{condi} and
\ref{eulerequ} is then

\begin{cor} \label{cor1}
Let $(M,g)$ be a compact Riemannian manifold of dimension $n \geq 3$ whose
Yamabe invariant is non-negative. We
assume that one of the following assumptions is true:
\begin{enumerate}[\theobull]
\item $\mu_1(M,g)>0$, $(M,g)$ is not locally conformally flat and $n \geq
11$;
\medskip
\item $\mu_1(M,g)=0$, $(M,g)$ is not locally conformally flat and $n \geq
9$.
\end{enumerate}

Then, there exists a nodal solution of Yamabe equation (\ref{eqyam}).
\end{cor}

{\it Acknowledgement\\}
The authors want to thank M. Ould Ahmedou for many 
interesting conversations
about nodal solutions of the Yamabe equation.
His large knowledge about such problems was a stimulating inspiration
for this article.
The author are also extremely obliged to Fr\'ed\'eric Robert for having
pointed out a little mistake in the first version of this paper.

\section{Variational characterization of $\mu_2(M,g)$} \label{varia}
\subsection{Notation}
In the whole article we will use the following notations
  $$L_+^N(M):=\left\{u\in L^N(M)\,|\,u\geq 0, \quad u\not \equiv 0\right\}.$$

\subsection{Grassmannians and the min-max principle}
Let $\Gr{k}{}{C^\infty(M)}$ be the $k$-dimensional \emph{Grassmannian} in 
$C^\infty(M)$, i.e.\ the set of all $k$-dimensional subspaces of $C^\infty(M)$.
The Grassmannian is an important ingredient in the min-max characterization of $\la_k(g)$  
  $$\la_k(L_{\ti{g}}):=\inf_{V\in \Gr{k}{}{C^\infty(M)}} \sup_{v\in \Vzero} {\int (L_{\ti g} v)v\, dv_{\ti g}\over \int_M v^2 \,dv_{\ti g}}.$$
We will also need a slightly modified Grassmannian. 
For any $u\in L_+^N(M)$ we define 
$\Gr{k}u{C^\infty}$ to be the set of all $k$-dimensional subspaces of 
$C^\infty(M)$, such that the restriction operator to $M\setminus u^{-1}(0)$
is injective. More explicitly, we have $\spann(v_1,\ldots,v_k) \in\Gr{k}{u}{C^\infty(M)}$ 
if and only if $v_1|_{M\setminus u^{-1}(0)},\ldots,v_k|_{M\setminus u^{-1}(0)}$
are linearly independent.
Sometimes it will be convenient to use the equivalent statement that 
the functions 
$u^{{N-2\over 2}}v_1,\ldots,u^{{N-2\over 2}}v_k$ are linearly independent.

Similarly, by replacing $C^\infty(M)$ by $H_1^2(M)$ we obtain the definitions
of $\Gr{k}{}{H_1^2(M)}$ and  $\Gr{k}u{H_1^2(M)}$.

\subsection{The functionals}
For all $u \in L_+^N(M)$, $v \in H_1^2(M)$ such that $u^{\frac{N-2}{2}} v
\not\equiv 0$, we set 
  $$F(u,v) =  \frac{\int_M c_n | \nabla v|^2 + S_g v^2 \,dv_g }{\int_M v^2 
    u^{N-2}\,dv_g} {\left( \int_M u^N \,dv_g \right)}^{\frac{2}{n}}.$$

\subsection{Variational characterization of $\mu_2(M,g)$}
The following characterization will be of central importance for our article.

\begin{prop}
We have 
\begin{eqnarray} \label{defmu}  
\mu_k(M,g)= \inf_{{\ss u\in L_+^N(M)\atop \ss V\in \Gr{k}u{H_1^2(M)}}}
\sup_{v\in\Vzero}  F(u, v) 
\end{eqnarray} 
\end{prop}

\proof{}

Let $u$ be  a smooth positive function on $M$. 
For all smooth functions $f$, 
$f \not\equiv 0$, we set $\tilde{g} = u^{N-2} g$ ($N= \frac{2n}{n-2}$)
and 
$$F'(u,f)= \frac{\int_M f L_{ \tilde{g}} f \,dv_{\tilde{g}}} {\int_M f^2  \,dv_{\tilde{g}}}.$$
The operator $L_g$ is conformally invariant (see \cite{hebey:97}) in the 
following sense: 
\begin{eqnarray} \label{conf_inv}
u^{N-1} L_{\tilde{g}} ( u^{-1} f) = L_g (f)
\end{eqnarray}
Together with the fact that 
\begin{eqnarray} \label{vol_elem}
dv_{\tilde{g}}= u^N \,dv_g,
\end{eqnarray}
 we get that 
$$F'(u,f)= \frac{\int_M (uf) L_{g}(uf)\,dv_g }{\int_M (u f)^2 u^{N-2}  
\,dv_g}  .$$

Using the min-max principle, we can write that 

$$\lambda_k (\tilde{g}) = \inf_{ V \in \Gr{k}u{H_1^2(M)} }\sup_{f \in \Vzero} F'(u,f)$$
Now, replacing $uf$ by $v$, we obtain that 
\begin{eqnarray} \label{def_lambda2}
\lambda_k (\tilde{g}) = \inf_{V\in \Gr{k}{}{H_1^2(M)}} \sup_{v \in \Vzero} \frac{\int_M v L_{g} v \,dv_g }{\int_M v^2 u^{N-2}  \,dv_g}.
\end{eqnarray} 

Using the definition of $\mu_2$ and 
$\Vol_{\tilde{g}}(M) = \int_M u^N \,dv_g$, we derive

$$\mu_k(M,g)= \inf_{\ss u\in L_+^N(M)\atop \ss V\in \Gr{k}u{C^\infty(M)}}
 \sup_{v\in\Vzero}  F(u, v) $$ 
The result follows immediately.

\section{Generalized metrics and the Euler-Lagrange equation}

\subsection{A regularity result}
We will need the following result.
\begin{lem} \label{regu}
Let $u \in L^N(M)$ and $v \in H_1^2(M)$. We assume that 
$$L_g v = u^{N-2} v$$
holds in the sense of distributions.
Then, $v \in L^{N+\ep}(M)$ for some $\ep >0$.
\end{lem}
This result is well known for the standard Yamabe equation. 
Proofs for the standard Yamabe equation 
can be found in \cite{trudinger:68} and \cite{hebey:97}, and the modifications
for proving Lemma~\ref{regu} are obvious. 
Unfortunately, \cite{trudinger:68} contains some typos, and  the book 
\cite{hebey:97} is difficult to obtain. This is why we included a proof in
the appendix for the convenience of the reader.

\subsection{The $k$-th eigenvalue of the Yamabe operator for a generalized
  metric} 
On a given Riemannian manifold $(M,g)$ we say that 
$\tilde{g} = u^{N-2} g$, $u\in L_+^N(M)$, is a 
\emph{generalized metric} conformal to $g$. 
For a  generalized metric $\tilde{g}$, we can define 
\begin{eqnarray}\label{def_lambdad}
\lambda_k (\tilde{g}) = \inf_{V\in \Gr{k}u{H_1^2(M)}} \sup_{v \in \Vzero} \frac{\int_M v L_{g} v \,dv_g }{\int_M v^2 u^{N-2}  \,dv_g}.
\end{eqnarray} 
%

\begin{prop} \label{la1la2} 
For any $u\in L_+^N$, $\ti g=u^{N-2}$ there exist two functions $v,w$ 
belonging to $H_1^2(M)$  
with $v \geq 0$ and such that in the sense of distributions.
\begin{eqnarray} \label{eqvl}
L_g v = \la_1 (\tilde{g})   u^{N-2} v
\end{eqnarray}
and
\begin{eqnarray} \label{eqwl}
L_g w = \la_2 (\tilde{g})   u^{N-2} w.
\end{eqnarray}
Moreover, we can normalize $v,w$ by 
\begin{eqnarray} \label{vwort}
\int_M u^{N-2} v^2 \,dv_g = \int_M u^{N-2} w^2 \,dv_g=1 \hbox{ and } 
\int_M u^{N-2} v w \,dv_g =0
\end{eqnarray}
\end{prop} 
For $k=2$ the  infimum in formula~(\ref{def_lambda2}) over all subspaces 
$V\in \Gr2{u}{H_1^2(M)}$ is attained by 
$V= \spann (v,w)$ and the
supremum over the functions in $V\setminus \{0\}$ is attained by $w$.
The reader should pay attention to the fact that the space $V$ is in
general non unique. As one can check, if $w$ changes the sign then 
the supremum over all $v\in V =\spann (v,w )\setminus\{0\}$ and the supremum
over all $v\in V_1= \spann ( w, |w| )\setminus\{0\}$ coincide.

{}From section  (\ref{varia}), we get   
$$\mu_2(M,g) = \inf_{\tilde{g} \in \ol{[g]}}   \lambda_2 (\tilde{g}).$$
Hence, $\mu_2(M,g)$ can be attained by a regular metric, or by 
a generalized metric or it can be not attained at all. 
These questions are discussed in
Section~\ref{properties}.  Let us now prove Proposition \ref{la1la2}. \\

 {\bf Proof of  Proposition \ref{la1la2}:} Let $(v_m)_m$ be a
minimizing sequence for $\la_1(\tilde{g}) $ i.e. a sequence $v_m\in H^2_1(M)$  
such that 
$$\lim_{m\to \infty}    
\frac{\int_M c_n |\nabla v_m|^2 + S_g v_m^2    \,dv_g}{\int_M  u^{N-2} v_m^2
  \,dv_g}= \la_1(\tilde{g}).$$
It is well known that $(|v_m|)_m$ is also a minimizing sequence. Hence, we can
assume that $v_m \geq 0$. 
If we normalize $v_m$ by $\int_M u^{N-2} v_m^2 \,dv_g=1$, then 
$(v_m)_m$ is bounded in $H^2_1(M)$ and after restriction to a subsequence 
we may assume
that there exists $v \in H^2_1(M)$, $v \geq 0$ 
 such that $v_m \to v$ weakly in
$H^2_1(M)$, strongly in $L^2(M)$ and almost everywhere. 
If $u$ is smooth, then 
\begin{eqnarray} \label{eq.limit}
\int_M u^{N-2} v^2 \,dv_g = \lim_ m \int_M u^{N-2} v_m^2 \,dv_g=1
\end{eqnarray}
and by standard arguments, $v$ is a non-negative minimizer of the
functional associated to $\la_1(\tilde{g})$.  
We must show that (\ref{eq.limit}) still holds if $u \in L_+^N(M) $. Let $A>0$
be a large real number and set $u_A = \inf(u,A)$. Then, using the H\"older
inequality, we write 
\begin{eqnarray*}
\left|  \int_M u^{N-2} \left(v_m^2 - v^2\right) \,dv_g \right| & 
\leq &  \left( \int_M u_A^{N-2} |v_m^2 - v^2| \,dv_g + \int_M (u^{N-2}
       - u_A^{N-2}) (|v_m|+|v|)^2 \,dv_g\right) \\
        & \leq & A \int_M |v_m^2 - v^2| \,dv_g\\ 
        &&{}+ {\left( \int_M (u^{N-2}-u_A^{N-2})^\frac{N}{N-2} \,dv_g
\right)}^{\frac{N-2}{N}}  {\left( \int_M (|v_m|+|v|)^N \,dv_g
  \right)}^{\frac{2}{N}}.
\end{eqnarray*}
By Lebesgue's theorem we see that 
  $$\lim_{A \to +\infty}  \int_M (u^{N-2}-u_A^{N-2})^\frac{N}{N-2} \,dv_g = 0. $$
Since $(v_m)_m$ is bounded in $H_1^2(M)$, it is bounded in $L^N(M)$ and
hence  there exists $C>0$ such that 
$\int_M (|v_m|+|v|)^N \,dv_g \leq C$. By strong convergence in $L^2(M)$, 
  $$\lim_m   \int_M |v_m^2 - v^2| \,dv_g =0.$$
Equation~\eref{eq.limit} easily follows and 
$v$ is a non-negative minimizer of the
functional associated to $\la_1(\tilde{g})$. Writing the Euler-Lagrange
equation of $v$, we find that $v$ satisfies equation (\ref{eqvl}).
Now, we define 
$$\la_2'(\tilde{g}) =  \inf \frac{\int_M  c_n |\nabla w|^2 + S_g w^2  \,dv_g}{\int_M  u^{N-2} |w|^2 \,dv_g}$$
where the infimum is taken over smooth functions $w$ such that 
$u^{\frac{N-2}{2} } w \not\equiv 0$ and such that 
$\int_M u^{N-2} vw \,dv_g= 0 $. With the same method, we find a minimizer
$w$ of this problem that satisfies (\ref{eqwl}) with $\la_2'(\tilde{g})$
instead of $\la_2(\tilde{g})$. However, it is not difficult to see that 
$\la_2'(\tilde{g})=\la_2(\tilde{g})$ and Proposition~\ref{la1la2} easily follows.

\subsection{Euler-Lagrange equation of a minimizer of $\la_2\Vol^{2/n}$}

\begin{lemma}\label{lem.EL}
Let $u\in L_+^N(M)$ with $\int u^N=1$. 
Suppose that $w_1,w_2\in H_1^2(M)\setminus\{0\}$, $w_1,w_2\geq 0$ satisfy
\begin{align}
 \int (c_n|\na w_1|^2 +\Scal_g w_1^2)\, dv_g & \leq \mu_2(M,g)\,\int u^{N-2}w_1^2 \label{ineq.v1}\\
 \int (c_n|\na w_2|^2 +\Scal_g w_2^2)\, dv_g & \leq \mu_2(M,g)\,\int u^{N-2}w_2^2\label{ineq.v2}
\end{align}
and suppose that $(M\setminus w_1^{-1}(0))\cap (M\setminus w_2^{-1}(0))$ 
has measure zero. Then $u$ is a linear combination of $w_1$ and~$w_2$
and we have equality in \eref{ineq.v1} and~\eref{ineq.v2}.
\end{lemma}

\proof{}
We let $\bar{u} = a w_1 + b w_2$ 
where $a,b>0$ are chosen such that
\begin{eqnarray} \label{v1v2} 
\frac{a^{N-2}}{b^{N-2}}\;
  \frac{\int_M u^{N-2} w_1^2 \,dv_g}{\int_M u^{N-2} w_2^2 \,dv_g}=
  \frac{\int_M w_1^N \,dv_g}{\int_M w_2^N\,dv_g}
\end{eqnarray}
and 
\begin{eqnarray} \label{baru}
\int_M \bar{u}^N \,dv_g=a^N\int_M w_1^N + b^N \int w_2^N=1.
\end{eqnarray} 
Because of the variational characterization of $\mu_2$
we have 
\begin{eqnarray} \label{mu<}
\mu_2(M,g) \leq \sup_{(\la, \mu) \in \mR^2 \setminus \{(0,0)\}}  
F(\bar{u}, \la w_1+ \mu w_2)
\end{eqnarray}
By \eref{ineq.v1},\eref{ineq.v2} and \eref{baru}, and
since  $(M\setminus w_1^{-1}(0))\cap (M\setminus w_2^{-1}(0))$ 
has measure zero 
\begin{eqnarray}
 F(\bar{u}, \la w_1+ \mu w_2) & = &
\frac{  \la^2 \int_M \left(c_n {|\nabla w_1|}^2 +S_g w_1^2\right) \,dv_g +  
\mu^2 \int_M\left( c_n {|\nabla
    w_2|}^2 +S_g w_2^2\right) \,dv_g }{\la^2  \int_M |\bar{u}|^{N-2} w_1^2 
    \,dv_g +\mu^2  \int_M |\bar{u}|^{N-2} w_2^2 \,dv_g}\nonumber \\
&\leq & \mu_2(M,g)  \frac{ \la^2\int_M u^{N-2} w_1^2 \,dv_g + \mu^2 \int_M
  u^{N-2} w_2^2 \,dv_g}{\la^2  a^{N-2} \int_M w_1^N \,dv_g + 
  \mu^2 b^{N-2} \int_M
  w_2^N \,dv_g}.\label{ineq.F}
\end{eqnarray}
As one can check, relation (\ref{v1v2}) implies that this expression does
not depend on $\la, 
\mu$. Hence, setting $\la=a $ and $\mu= b$, the denominator 
is $1$, and we get 
\begin{eqnarray*}
\sup_{ (\la, \mu) \in \mR^2 \setminus \{(0,0)\}} 
F(\bar{u}, \la w_1+ \mu w_2) & \leq &  \mu_2(M,g) \int_M u^{N-2}(a^2 w_1^2 + 
b^2 w_2^2) \,dv_g \\
& = & \mu_2(M,g) \int_M u^{N-2} \bar{u}^2 \,dv_g.
\end{eqnarray*}
By H\"older inequality, 
\begin{eqnarray}
\sup_{ (\la, \mu) \in \mR^2 \setminus \{(0,0)\}}  
F(\bar{u}, \la w_1+ \mu w_2) \leq  \mu_2(M,g) {\left(  \int_M u^N \,dv_g
   \right)}^{\frac{N-2}{N}} {\left(  \int_M {\bar{u}}^N \,dv_g
   \right)}^{\frac{2}{N}} =\mu_2(M,g).\label{ineq.hoelder}
\end{eqnarray}
Inequality  (\ref{mu<}) implies that we have both equality in the 
H\"older inequality of \eref{ineq.hoelder} and in~\eref{ineq.F}.
The equality in the H\"older inequality implies that there exists a constant 
$c>0$ such that $u= c \bar{u}$ almost everywhere. Moreover, since 
$\int u^N = \int \bar{u}^N =1$, we have 
$u = \bar{u}= a w_1 + b w_2 $.
The equality in \eref{ineq.F} implies inequality in 
\eref{ineq.v1} and~\eref{ineq.v2}. 
\qed

\begin{theorem}[Euler-Lagrange equation] \label{theo.limit}
Assume that $\mu_2(M,g)\neq 0$ and that 
$\mu_2(M,g)$ is attained by a generalized metric $\tilde{g} =
u^{N-2} g$ with $ u \in L_+^N(M)$. 
Let $v,w$ be as in Proposition~\ref{la1la2}. 
Then, $u = |w|$. In particular,
\begin{eqnarray} \label{eqlim}
L_g w= \mu_2(M,g) |w|^{N-2} w
\end{eqnarray}
Moreover, $w$ has alternating sign and $w \in
C^{3,\alpha}(M)$ ($\alpha \leq N-2$).
\end{theorem}

\begin{rem}
Assume that $\mu_2(M,g)$ is equal to $0$ and is attained by a generalized metric $g'$,
then, using the conformal invariance of the Yamabe operator, 
it is easy to
check that for all generalized metrics $\tilde{g}$ conformal to $g'$, we have 
$\lambda_2(\tilde{g}) = 0$. Consequently, each metric conformal to $g$ is a
minimizer for $\mu_2(M,g)$ and Theorem \ref{theo.limit} is
always false in this case. However, we will still get a nodal solution
of~\eref{eqyam} if
$\mu_2(M,g)=0$. Indeed, by Theorem~\ref{attain} and the remark above, 
$\la_2(g)=0$. Let $w$ be an eigenfunction associated
to $\la_2(g)$. We have $L_g w = 0$. Then, we have a solution
of  \eref{eqlim}.
\end{rem} 

\begin{rem} 
Assume that $\mu_2(M,g) \not= 0$ and that $\mu_2(M,g)$ is attained by a
generalized metric. Let $w$ be the solution of equation (\ref{eqlim})
given by
Theorem \ref{theo.limit}. We let $\Om_+= \{x \in M \hbox{ s.t. } w(x)>0 \}$
and $\Om_-= \{x \in M \hbox{ s.t. } w(x)-0 \}$. Then, a immediate
consequence of Lemma \ref{lem.EL} is that 
$\Om_+$ and $\Om_-$ have exactly one connex component. 
\end{rem}

\proof{}
Without loss of generality, we can assume that $\int_M u^N \,dv_g =1$. 
By assumption we have 
$ \la_2(\tilde{g})= \mu_2(M,g)$. 
Let $v,w \in H_1^2(M)$ be some functions satisfying equations
(\ref{eqvl}), (\ref{eqwl}) and relation (\ref{vwort}).

\begin{step} \label{ste1}
We have $\la_1(\tilde{g}) < \la_2
(\tilde{g})$. 
\end{step}
We assume that $\la_1(\tilde{g}) = \la_2
(\tilde{g})$. Then, after possibly replacing  
$w$ by a linear combination of $v$
and $w$, we can assume that the function 
$u^{\frac{N-2}{2}} w$ changes the sign. 
We apply Lemma~\ref{lem.EL}  for $w_1 := \sup(w,0)$
and   $w_2:= \sup(-w,0)$.
We obtain the existence of $a,b>0$ with
$u = a w_1 + b w_2$. Now, by Lemma~\ref{regu},
$w \in L^{N+\ep}(M)$. By a standard bootstrap argument, equation (\ref{eqwl})
shows that  
$w \in C^{2, \al}(M)$ for all $\al \in ]0,1[$. It follows that $u\in C^{0,
  \al}(M)$ for all $\al \in ]0,1[$. 
Now, since $\la_1(\tilde{g}) = \la_2
(\tilde{g})$ and by definition of $\la_1(\tilde{g})$, $w$ is a minimizer of
the functional $\bar{w} \mapsto F(u,\bar{w})$ among the functions belonging to 
$H_1^2(M)$ and such that $u^{\frac{N-2}{2}} \bar{w} \not\equiv 0$. Since 
$F(u,w) = F(u,|w|)$, we see that $|w|$ is a minimizer for the functional
associated to $\la_1(\tilde{g})$ and hence, writing the Euler-Lagrange
equation of the problem, $w$ 
satisfies the same equation 
as $w$. As a consequence,  $|w|$ is
$C^2(M)$. By the maximum principle, we get $|w| >0$ everywhere. This
is false. Hence, the step is proved.\\

\begin{step}
The function~$w$ changes the sign.
\end{step}

Assume that $w$ does not change the sign, i.e.~ after possibly replacing $w$
by $-w$, we have $w\geq 0$. Using~\eref{vwort} we see that 
$(M\setminus v^{-1}(0))\cap (M\setminus w^{-1}(0))$ has measure zero.
Setting $w_1:=v$ and $w_2:=w$ we have \eref{ineq.v1} and~\eref{ineq.v2}.
While we have equality in \eref{ineq.v2}, Step~1 implies that 
inequality \eref{ineq.v1} is strict.
However using Lemma~\ref{lem.EL} we can derive equality in \eref{ineq.v1}.
Hence we obtain a contradiction, and the step is proved.


\begin{step} 
There exists $a,b >0$ such that $u=  a \sup(w,0) + b \sup(-w,0)$.
Moreover, $w\in C^{2,\al}(M)$ and   
$u\in C^{0,\al}(M)$ for all $\al \in ]0,1[$.
\end{step}

As in the proof of Step~1 we apply Lemma~\ref{lem.EL} 
for $w_1:=\sup(w,0)$ and $w_2:=\sup(-w,0)$. 
We obtain the existence of $a,b >0$ such that $u=  a w_1 + b w_2$.
As in Step~1 we get that $w\in C^{2,\al}(M)$ and   $u\in C^{0,\al}(M)$ 
for all $\al \in ]0,1[$. This proves the present step.

\begin{step}
Conclusion.
\end{step}
 Let $h \in C^{\infty}(M)$ whose support is contained in $M \setminus
 \{u^{-1}(0) \}$.    
For $t$ close to $0$, set $u_t = |u + th|$. Since $u>0$ on the support of
$h$ and since $u$ is continuous (see last step), we have  for $t$ close to
$0$,
$u_t=u+th$. 
As $\spann(v,w)\in \Gr{2}u{H_1^2(M)}$ we obtain using \eref{defmu} 
for all $t$
  $$\mu_2(M,g) \leq  \sup_{(\lambda, \mu) \in \mR^2 \setminus \{ (0,0) \}} 
    F(u_t,\la v + \mu w).$$
Equations (\ref{eqvl}), (\ref{eqwl}), and relation (\ref{vwort})
yield 
\begin{eqnarray*}
F(u_t,\la v + \mu w
 )& = & \frac{\la^2 \la_1(\ti g) \int_M u^{N-2} v^2 \,dv_g 
  + \mu^2 \la_2(\ti g) \int_M u^{N-2} w^2
  \,dv_g}
  {\la^2 \int_M u_t^{N-2} v^2 \,dv_g +  2 \la
  \mu  \int_M u_t^{N-2}  v w
  \,dv_g + \mu^2 \int_M u_t^{N-2} w^2 \,dv_g } 
  { \left( \int_M u_t^N \,dv_g \right)}^{\frac{2}{n}}\\
  & =&  \frac{\la^2 \la_1(\ti g)  + \mu^2 \la_2(\ti g)}{\la^2 a_t 
        + \la\mu  b_t + \mu^2 c_t} 
        {\left( \int_M |u_t|^N \,dv_g \right)}^{\frac{2}{n}}.
\end{eqnarray*} 
where
$$a_t= \int_M u_t^{N-2} v^2 \,dv_g,$$
$$b_t = 2  \int u_t^{N-2}  v w
  \,dv_g$$
and 
$$c_t = \int_M u_t^{N-2} w^2 \,dv_g.$$
The functions $a_t$, $b_t$ and $c_t$ are 
smooth for $t$ close to $0$, 
furthermore $a_0 =  c_0=1$ and $b_0=0$.
The function $f(t,\al):= F(u_t,\sin(\al) v + \cos(\al) w)$ is smooth for 
small $t$.
Using $\la_1(\ti g)<\la_2(\ti g)$
one calculates
\begin{align*}
{\pa\over \pa\al}\,f(0,\al)    &=0 &\Leftrightarrow\qquad&\al\in {\pi\over 2}\mZ\\
{\pa^2\over \pa\al^2}\,f(0,\al)&<0 &\mbox{for}     \qquad&\al\in \pi\mZ\\ 
{\pa^2\over \pa\al^2}\,f(0,\al)&>0 &\mbox{for}     \qquad&\al\in \pi\mZ+{\pi\over 2}
\end{align*} 
Applying the implicit function theorem to ${\pa f \over \pa\al}$ at the
point $(0,0)$,  we see that  
there is a smooth function $t\mapsto\alpha(t)$, 
defined on a neighborhood of $0$
with $\al(0)=0$ and
   $$f(t,\al(t))=\sup_{\al\in \mR}f(t,\al)= 
     \sup_{(\lambda, \mu) \in \mR^2 \setminus \{ (0,0) \}} 
      F(u_t,\la v + \mu w).$$ 
As a consequence 
  $${d\over dt}|_{t=0}\sin^2\al(t)={d\over dt}|_{t=0}\cos^2\al(t)={d\over dt}|_{t=0}(\sin^2\al(t)a_t)={d\over dt}|_{t=0}(\sin\al(t)\cos\al(t)b_t)=0.$$
Hence, ${d\over dt}|_{t=0}\,f(t,\al(t))$ exists and we have 
\begin{eqnarray*} 
{d\over dt}|_{t=0}\,f(t,\al(t)) & = & \la_2(M,\ti g)\left(- \frac{d}{dt}|_{t=0} c_t + \frac{d}{dt}|_{t=0}
{\left( \int_M |u_t|^N \,dv_g \right)}^{\frac{2}{n}} \right) \\
& = & \la_2(M,\ti g) (N-2)\left( - \int_M u^{N-3} h w^2 \,dv_g + \int_M u^{N-1} h
  \,dv_g \right).
\end{eqnarray*}
By definition of $\mu_2(M,g)$, $f$ admits a minimum in $t= 0$. As 
$\la_2(M,\ti g)=\mu_2(M,g)\neq 0$ we obtain
$$ \int_M u^{N-3} h w^2 \,dv_g = \int_M u^{N-1} h
  \,dv_g.$$
Since $h$ is arbitrary (we just have to ensure that its support is contained
in $M \setminus \{ u^{-1}(0)\}$), we get that 
$u^{N-3}w^2 =  u^{N-1}$ on $M \setminus \{ u^{-1}(0)\}$,
hence $u=|w|$ on $M \setminus \{ u^{-1}(0)\}$.
Together with Step~3, we get $u=|w|$ everywhere. 
This proves theorem~\ref{theo.limit}.\qed

\section{A sharp Sobolev inequality related to $\mu_2(M,g)$}

\subsection{Statement of the results}

For any compact Riemannian manifold $(M,g)$ of dimension $n \geq 3$, Hebey and Vaugon have shown in (\cite{hebey.vaugon:96}) that there exists $B_0(M,g) >  0$ 
such that 
  $$ \mu_1(\mS^n)= n(n-1)\,\om_n^{\frac{2}{n}}  
     = \inf_{u \in H_1^2(M)\setminus \{0\}}  \frac{
      \int_M c_n |\nabla u |^2 + B_0 \int_M u^2 \,dv_g}
     {{\left(\int_M u^N \,dv_g \right)}^{\frac{2}{n}}} \eqno{\rm (S)}$$
where $\om_n$ stands for the volume of the standard $n$-dimensional sphere $\mS^n$ and where $\mu_1(\mS^n)$ is the Yamabe invariant of $\mS^n$. 

This inequality is strongly related to the resolution of the Yamabe problem. 
It allows to avoid concentration for the minimizing sequence of $\mu_1 (M,g)$. 
For the minimization of $\mu_2(M,g)$, this inequality 
is not sufficient and another one must be constructed. The following result is adapted to the problem of minimizing $\mu_2(M,g)$. 

\begin{theorem} \label{sobolev}
On a  compact connected 
Riemannian manifold $(M,g)$ of dimension $n \geq 3$
we have 
$$ 2^{\frac{2}{n}} \mu_1(\mS^n)  = \inf_{\ss u  \in L_+^N(M)\atop\ss V\in \Gr{k}u{H_1^2(M)}} \sup_{v\in \Vzero}
 \frac{\left(\int_M c_n |\nabla v  |^2 \,dv_g +B_0(M,g) 
\int_M  v^2 \,dv_g \right) 
{\left(  \int_M u^N \,dv_g \right)}^{\frac{2}{N}}}{ 
\int_M u^{N-2} v^2 \,dv_g} \eqno{\rm (S_1)}$$
where $B_0(M,g)$ is given by inequality {\rm (S)}. 
\end{theorem}

We present now two  corollaries  of Theorem~\ref{sobolev}.

\begin{cor} \label{muS}
For the standard $n$-dimensional sphere we have 
$ \mu_2 (\mS^n)  =    2^{\frac{2}{n}} \mu_1(\mS^n)$.   
\end{cor}

\begin{cor} \label{muR}
For all $u \in C^{\infty}_c(\mR^n)$ and $V\in \Gr{2}u{C^\infty_c(\mR^n)}$
we have
$$ 2^{2/n} \mu_1(\mS^n)  \leq \sup_{v\in \Vzero} 
 \frac{\left(\int_{\mR^n} c_n |\nabla v  |^2 \,dv_g \right) {\left(  \int_{\mR^n} |u|^N \,dv_g \right)}^{\frac{2}{N}}}{ \int_{\mR^n} |u|^{N-2} v^2 \,dv_g} $$
\end{cor}

\subsection{Proof of theorem \ref{sobolev}}
The functional 
  $$G(u,v):=\frac{\left(\int_M c_n |\nabla  v |^2 \,dv_g +B_0(M,g) 
    \int_M v^2 \,dv_g \right) {\left(  \int_M u^N \,dv_g \right)}^{\frac{2}{N}}}
    { \int_M u^{N-2} v^2 \,dv_g}$$
is continuous on $L_+^N(M)\times (H_1^2(M)\setminus\{0\})$. 
As a consequence $I(u,V):= \sup_{v\in V \setminus \{0\}} G(u,v)$ depends continuously
on $u\in L_+^N(M)$ and $V\in \Gr{2}u{H_1^2(M)}$.
Thus, in order to show the theorem it is sufficient to show that
$I(u,V)\geq 2^{2/n} \mu_1(\mS^n)$ for all smooth $u>0$ and 
$V\in \Gr{2}{}{C^\infty(M)}$.
Without loss of 
generality,  we can assume
\begin{eqnarray} \label{u=1}
\int_M u^N \,dv_g =1.
\end{eqnarray} 
The operator 
$v\mapsto P(v):= c_n u^{2-N\over 2}\Delta (u^{2-N\over 2}v) +  B_0(M,g)u^{2-N}v$ is an 
elliptic operator on $M$, and $P$ is self-adjoint with respect to the $L^2$-scalar
 product. Hence, $P$ has discrete spectrum 
$\la_1\leq \la_2\leq \ldots$ and the corresponding eigenfunctions 
$\phi_1,\phi_2,\ldots$ are smooth.
Setting $v_i:= u^{2-N\over 2}\phi_i$ we obtain
  $$\left(c_n\Delta+B_0\right)(v_i)= \la_i u^{N-2}v_i$$
  $$\int u^{N-2}v_iv_j\,dv_g=0\qquad\mbox{if $\la_i\neq \la_j$}.$$
The maximum principle implies that an eigenfunction to the smallest eigenvalue
$\la_1$ has no zeroes. Hence $\la_1<\la_2$, and we can assume $v_1>0$.

We define $w_+:=a_+\sup(0,v_2)$ and $w_-:=a_-\sup(0,-v_2)$, where we choose
$a_+,a_->0$ such that

$$\int_M u^{N-2} w_-^2 \,dv_g = \int_M u^{N-2} w_+^2 \,dv_g =1.$$

We let $\Om_-= \{ w <0 \} $ and $\Om_+= \{ w \geq 0 \} $. 
By H\"older inequality,
 
\[ \begin{array}{ccl}
2 & = & \int_M u^{N-2} w_-^2 \,dv_g + \int_M u^{N-2} w_+^2 \,dv_g \\
& \leq &
{\left(\int_{\Om_-} u^N \,dv_g \right)}^{\frac{N-2}{N}}  {\left(\int_M w_-^N 
\,dv_g \right)}^{\frac{2}{N}}
+  {\left(\int_{\Om_+} 
u^N \,dv_g \right)}^{\frac{N-2}{N}}  {\left(\int_M w_+^N 
\,dv_g \right)}^{\frac{2}{N}}. \end{array} \]
Using the sharp Sobolev inequality $(S)$, we get that 
\begin{eqnarray} \label{usingS}
2 \mu_1(\mS^n) &\leq & {\left(\int_{\Om_-} u^N \,dv_g \right)}^{\frac{N-2}{N}} \int_M w_- u^{N-2\over 2}P \left( u^{N-2\over 2}\, w_-\right) \,dv_g\\
& + & {\left(\int_{\Om_+} u^N \,dv_g \right)}^{\frac{N-2}{N}} \int_M w_+ u^{N-2\over 2}P \left( u^{N-2\over 2}\, w_+\right) \,dv_g
\end{eqnarray}
Since $w_-$ resp.\ $w_+$ are some multiples of $w$ on $\Om_-$ resp.\ $\Om_+$, 
they satisfy 
the  same equation as $w$. Hence, we get that 
\[ \begin{array}{ccl}
2 & = &  \mu_1(\mS^n)^{-1}  \la_2 \left(
{\left(\int_{\Om_-} u^N \,dv_g \right)}^{\frac{N-2}{N}} 
\int_M u^{N-2} w_-^2 \,dv_g 
+ {\left(\int_{\Om_+} u^N \,dv_g \right)}^{\frac{N-2}{N}} 
\int_M u^{N-2} w_+^2 \,dv_g \right) \\
& = &  \mu_1(\mS^n)^{-1}  \la_2 \left( {\left(\int_{\Om_-} u^N \,dv_g \right)}^{\frac{N-2}{N}} + {\left(\int_{\Om_+} u^N \,dv_g \right)}^{\frac{N-2}{N}} \right)
. \end{array} \]
Now, for any real non-negative numbers $a,b \geq 0$, the H\"older inequality yields
$$a + b \leq 
2^{\frac{2}{N}}{\left( a^{\frac{N}{N-2}} + b^ {\frac{N}{N-2}} 
\right)}^{\frac{N-2}{N}}  $$
We apply this inequality with 
$a = {\left(\int_{\Om_-} u^N \,dv_g \right)}^{\frac{N-2}{N}}$ and  
$b= {\left( \int_{\Om_+} u^N \,dv_g \right)}^{\frac{N-2}{N}} $. 
Using (\ref{u=1}), 
we obtain
$$ 2 \leq 2^{\frac{2}{N}} \mu_1(\mS^n)^{-1}  \la_2  \left( \int_{\Om_-} u^N \,dv_g + 
\int_{\Om_+} u^N \,dv_g \right)= 2^{\frac{2}{N}} \mu_1(\mS^n)^{-1}  \la_2. $$
We obtain 
$\la_2 \geq  2^{\frac{2}{n }} \mu(\mS^n)$. 
Since $\la_2 = I(u,\spann(v_1,v_2))$, this ends the proof of Theorem~\ref{sobolev}.

\subsection{Proof of Corollaries~\ref{muS} and \ref{muR}}

It is well known that $B_0(\mS^n)$ equals to the scalar curvature 
of~$\mS^n$, i.e.\ $B_0(\mS^n)=n(n-1)$.
Replacing $B_0(\mS^n)$ by its value and taking the infimum over $u,V$, 
the right hand term of inequality $(S_1)$ is exactly the variational
characterization of $\mu_2(\mS^n)$ (see equation \eref{defmu}). 
This proves that 
$\mu_2(\mS^n) \geq  2^{2/n} \mu_1(\mS^n)$. Corollary~\ref{muS} then 
follows from Theorem~\ref{upbound}. 
Since $\mR^n$ is conformal to $\mS^n \setminus \{p\}$ ($p$ is any point of $\mS^n$), we can use the conformal invariance to prove 
Corollary~\ref{muR}.

\section{Some properties of $\mu_2(M,g)$} \label{properties}
\subsection{Is $\mu_2(M,g)$ attained?} \label{attaine}

Let $(M,g)$ be an $n$-dimensional compact Riemannian manifold. 
The Yamabe problem shows that $\mu_1(M,g)$ is attained by a metric $\tilde{g}$ conformal to $g$. 
Some questions arise naturally concerning 
 $\mu_2(M,g)$:

{\bf 1-} Is  $\mu_2(M,g)$  
attained by a metric?

{\bf 2-} Is it possible that $\mu_2(M,g)$ is attained by a 
generalized metric?

In this section, we give answers to these questions.
The first result we prove is the following:

\begin{prop} \label{sUs}
Let $\mS^n \dot\cup \mS^n$ be the disjoint union of 
two copies of the sphere equipped with their standard metric. 
Then, $\mu_2(\mS^n \dot\cup \mS^n)= 2^{2/n} \mu_1(\mS^n)$ and 
it is attained by the canonical metric.
\end{prop}

\proof{}One computes 
  $$\la_2(\mS^n \dot\cup \mS^n)\,\Vol(\mS^n \dot\cup \mS^n)^{2/n}
    = 2^{2/n}\la_1(\mS^n)\,\Vol(\mS^n)^{2/n}=2^{2/n} \mu_1(\mS^n).$$
Hence $\mu_2(\mS^n \dot\cup \mS^n)\leq 2^{2/n} \mu_1(\mS^n)$ follows.

Now, let $\ti g$ be an arbitrary smooth metric on  $S^n \dot\cup S^n$.
We write $S^n_1$ for the first $S^n$ and $S^n_2$ for the second $S^n$.
Then $\la_2(S^n \dot\cup S^n,\ti g)$ is the minimum
of $\la_2(S^n_1,\ti g)$, $\la_2(S^n_2,\ti g)$ and 
$\max\{\la_1(S^n_1,\ti g),\la_1(S^n_2,\ti g) \}$.

It follows from Corollary~\ref{muS} that 
  $$ \la_2(S^n_1,\ti g)\,\Vol(S^n\dot\cup S^n,\ti g)^{2/n}\geq 
     \la_2(S^n_1,\ti g)\, \Vol(S^n_1,\ti g)^{2/n}\geq 
     2^{2/n }\mu_1(\mS^n),$$
and obviously we have the same for $\la_2(S^n_2,\ti g)$.

Summing 
  $$\la_1(S^n_i,\ti g)^{n/2}\geq \mu_1(\mS^n)^{n/2}\Vol(S^n_i,\ti g)$$
over $i\in\{1,2\}$, we obtain the remaining inequality
  $$\max\{\la_1(S^n_1,\ti g),\la_1(S^n_2,\ti g) \}\,
    \Vol(S^n\dot\cup S^n,\ti g)^{2/n}\geq  2^{2/n}  \mu_1(\mS^n),$$
and the proposition is proved.
\qed

Question~1 is solved by 
the following result.

\begin{prop} \label{notattained}
If $M$ is connected, then $\mu_2(M,g)$ cannot be attained by a metric.
\end{prop}
Indeed, otherwise by Theorem~\ref{theo.limit}, we would have that 
$u=|w|$ and hence $u$ cannot be positive. 
Theorem~\ref{attain} and the following result answer Question~2. 

\begin{prop} \label{notattainedS}
The invariant $\mu_2(\mS^n)$ is not attained by a generalized metric.
\end{prop}

This proposition immediately follows from Proposition~\ref{lowbound}.

\subsection{Some bounds of $\mu_2(M,g)$}
At first, we give an upper bound for $\mu_2(M,g)$.
\begin{theorem} \label{upbound} 
Let $(M,g)$ be an $n$-dimensional compact  Riemannian manifold with 
$\mu_1(M,g) \geq 00$. Then,
\begin{equation}\label{eq.upbound}
 \mu_2(M,g) \leq {( \mu_1(M,g)^{\frac{n}{2}}+
  \mu_1(\mS^n)^{\frac{n}{2}})}^{\frac{2}{n}}.
\end{equation}
This inequality is strict in the following cases:
\begin{enumerate}[\theobull]
\item $\mu_1(M,g) > 0$, $(M,g)$ is not locally conformally flat and $n \geq 11$;
\item $\mu_1(M,g)=0$, $(M,g)$ is not locally conformally flat and $n \geq 9$.
\end{enumerate}
\end{theorem}

\noindent {}From the solution of the Yamabe problem by Aubin and Schoen \cite{aubin:76,schoen:84}
we know that if $(M,g)$ is not conformally equivalent to $\mS^n$,
then $\mu_1(M,g)<\mu_1(\mS^N)$. Hence, \eref{eq.upbound} implies  the
following corollary.

\begin{cor}
Let $(M,g)$ be an $n$-dimensional compact connected Riemannian manifold
whose Yamabe invariant is non-negative. Then
$\mu_2(M,g) \leq \mu_2(\mS^n)$ with inequality if and only if $(M,g)$ is
conformally diffeomorphic to the sphere $\mS^n$. 
\end{cor}

These inequalities are very important, because they can be used
to avoid concentration of minimizing 
sequences for $\mu_2$, in a way which is similar to the resolution
of the Yamabe problem.

The following proposition gives a lower bound for~$\mu_2$.

\begin{prop} \label{lowbound}
Let $(M,g)$ be a $n$-dimensional compact Riemannian manifold 
whose  Yamabe invariant is non-negative. Then, 
 \begin{eqnarray} \label{lowerbound}
\mu_2(M,g)\geq 2^{2 \over n}  \mu_1(M,g).
\end{eqnarray}
Moreover, if $M$ is connected and if $\mu_2(M,g)$ is attained by a
generalized metric, then this inequality is strict.
\end{prop}

When $\mu_1(M,g) = 0 $, inequality (\ref{lowerbound}) is trivial. 
If $\mu_1(M,g) >0$, by a
possible chande of metric in the conformal class, we can assume that the
scalar curvature is positive.  The proof of inequality (\ref{lowerbound}) 
is exactly the
same as the one of Theorem~\ref{sobolev}. We just have to replace $B_0(M,g)$
by $S_g$. 
Moreover, if $M$ were connected and if $\mu_2(M,g)$ were attained by a
generalized metric, then inequality (\ref{usingS}) would be an equality and
we would have that $w_+$ or $w_-$ is a function for which equality in the 
Sobolev inequality $(S)$ is attained. 
By the maximum principle, we would get that $w_+$ or $w_-$
is positive on $M$ which is impossible.

\subsubsection{Proof of theorem \ref{upbound}}


\begin{lem} \label{estim}
For any $\al>2$, there is a  $C >0$ such that 
$$|a+b|^{\al} \leq a^{\al} + b^{\al} + C ( a^{\al -1}b  +  a b^{\al -1})$$
for all $a,b >0$.
\end{lem}

 {\bf Proof of Lemma \ref{estim}.}
Without loss of generality, we can assume that $a=1$. Then we set for $x > 0$,
$$f(x) = \frac{|1+x|^{\al}  - (1 + x^{\al})}{x^{\al-1} + x}.$$
One checks that $\lim_{x \to 0} f(x) = \lim_{x \to +\infty } f(x)=
\al$. Since $f$ is continuous, $f$ is bounded by a constant $C$ on
$\mR_+$. Clearly, this constant is the desired $C$ in  
inequality of Lemma \ref{estim}.\\

 {\bf Proof of Theorem \ref{upbound}.}
For $u \in H_1^2(M)\setminus \{0\}$ let
$$Y(u) = \frac{ \int_M c_n   |\nabla u|^2 +S_g u^2 \,dv_g}{{\left( \int_M
      |u|^N \,dv_g \right)}^{\frac{2}{N}}}$$
be the Yamabe functional of $M$.
The solution of the Yamabe
problem provides the existence of a smooth positive minimizer $v$ 
of~$Y$, and we can assume
\begin{eqnarray} \label{intv}
\int_M v^N \,dv_g = 1.
\end{eqnarray}
Then, $v$ satisfies the Yamabe equation
\begin{eqnarray} \label{yamequ}
L_g v = \mu_1(M,g) v^{N-1}.
\end{eqnarray}
Let $x_0 \in M$ be fixed and choose a system $(x_1,\cdots, x_n)$ of normal
coordinates at $x_0$. We note $r=dist_g(x_0,.)$. If $\delta >0$ is a small
fixed number, let $\eta$ be a smooth cut-off function such that
$0 \leq \eta \leq 1$, $\eta(B(x_0, \delta) ) = \{ 1 \}$ and $\eta( M
\setminus B(x_0, 2 \delta)= \{ 0 \}$, $|\nabla \eta|\leq 2/\delta$. 
Then, we can define for all $\ep >0$
$$v_{\ep} = C_{\ep} \eta (\ep + r^2)^{\frac{2-n}{n}}.$$
where  $C_{\ep}>0$ is such that 
\begin{eqnarray} \label{intve}
\int_M v_{\ep}^N \,dv_g = 1.
\end{eqnarray}
By standard computations (see \cite{aubin:76}) 
\begin{eqnarray} \label{testsph}
\lim_{\ep \to 0} Y(v_{\ep}) = \mu_1(\mS^n).
\end{eqnarray}
If $(M,g)$ is not locally conformally flat, if $g$ is well chosen in the
 conformal class and if $x_0$ is well chosen in $M$, it was also proven in
 \cite{aubin:76} that  there exists a constant $C(M)>0$ such that 
\begin{eqnarray} \label{test}
   Y(v_{\ep})=\left|\;  
   \begin{matrix}
     \mu_1(\mS^n) - C(M) \ep^2 + o(\ep^2)\hfill & \hbox{ if } n > 6 \\
     \mu_1(\mS^n) - C(M) \ep^2 |\ln(\ep)| + o(\ep^2 |\ln(\ep)|)\hfill 
     & \hbox{ if } n = 6.
   \end{matrix}
   \right.
\end{eqnarray} 
Moreover, it follows from \cite{aubin:76}
that $$ a \ep^{\frac{n-2}{4} } \leq C_{\ep} \leq b \ep^{\frac{n-2}{4} }$$
where $a,b >0$  are independent of $\ep$. If $p \geq 1$, standard computations
made in \cite{aubin:76} show 
that there exist  some constants $c,C>0$ independent of $\ep$ such that

\begin{eqnarray} \label{normp}
c \al_{p,\ep} \leq \int_M v_{\ep}^p \,dv_g \leq C \al_{p,\ep }
\end{eqnarray}
where 
\[ \al_{p,\ep} = \left| \begin{array}{lll} 
\ep^{\frac{2n - (n-2) p}{4}} & \hbox{if} & p> \frac{n}{n-2};\\
|\ln(\ep)| \ep^{\frac{n}{4}} & \hbox{if} & p= \frac{n}{n-2};\\
\ep^{ \frac{(n-2) p}{4}} & \hbox{if}& p< \frac{n}{n-2}
\end{array} \right. \]
Since the large inequality if easier to obtain, we only prove strict
inequality. Assume first that $\mu_1(M,g)>0$, that $(M,g)$ is not locally
conformally flat and that $n \geq 11$. We set, 
$$u_{\ep} = Y(v_{\ep})^{\frac{1}{N-2}}  v_{\ep} + \mu_1(M,g)^{\frac{1}{N-2}}  v.$$
Let us derive estimates for $F\big(u_{\ep},\la v_{\ep}+\mu v)\big)$. 
Let $(\lambda,\mu) \in \mR^2 \setminus \{(0,0)\}$.
Using (\ref{intv}), (\ref{intve}) and  the equation (\ref{yamequ}) of
$v$, we get that  
\begin{eqnarray*}
F(u_{\ep}, \la v_{\ep} + \mu v) & = &
 \frac{  \la^2 \int_M v_{\ep} L_g(v_{\ep})\,dv_g +  
\mu^2 \int_M v L_g(v) \,dv_g + 2 \la \mu  \int_M v_{\ep} L_g v \,dv_g }
{\la^2  \int_M |u_{\ep}|^{N-2}(\la v_{\ep} + \mu v)^2 \,dv_g}
{\left( \int_M u_{\ep}^N  \,dv_g \right)}^{\frac{2}{n}}.  
\end{eqnarray*}
\begin{eqnarray} \label{testfunc}
= \frac{  \la^2  Y(v_{\ep}) + \mu^2 \mu_1(M,g) + 2 \la \mu 
 \mu_1(M,g)  \int_M |v|^{N-2} v v_{\ep} \,dv_g}{ \la^2 \int_M |u_{\ep}|^{N-2}
v_{\ep}^2 \,dv_g  + \mu^2 \int_M |u_{\ep}|^{N-2} v^2 \,dv_g  + 2\la \mu 
\int_M |u_{\ep}|^{N-2} v v_{\ep} \,dv_g} {\left( \int_M u_{\ep}^N  \,dv_g
\right)}^{\frac{2}{n}}. 
\end{eqnarray}

Using the definition of $u_{\ep} $
\begin{eqnarray*}   
\la^2 \int_M |u_{\ep}|^{N-2}
v_{\ep}^2 \,dv_g  + \mu^2 \int_M |u_{\ep}|^{N-2} v^2 \,dv_g  + 2\la \mu 
\int_M |u_{\ep}|^{N-2} v v_{\ep} \,dv_g 
\end{eqnarray*}
\begin{eqnarray*} 
\;&\ \geq   \la^2 Y(v_{\ep}) \int_M |v_{\ep}|^{N} \,dv_g + \mu^2 
 \mu_1(M,g)  \int_M |v|^{N} \,dv_g +  2\la \mu 
\int_M |u_{\ep}|^{N-2} v v_{\ep} \,dv_g \\
\; & = \la^2 Y(v_{\ep})+ \mu^2 
 \mu_1(M,g) +  2\la \mu 
\int_M |u_{\ep}|^{N-2} v v_{\ep} \,dv_g.
\end{eqnarray*}
If $\la \mu \geq 0$, we have 
$$ 2\la \mu 
\int_M |u_{\ep}|^{N-2} v v_{\ep} \,dv_g \geq 2\la \mu \mu_1(M,g) 
\int_M v^{N-2}  v_{\ep} \,dv_g.$$
This implies that 
$$ \frac{  \la^2  Y(v_{\ep}) + \mu^2 \mu_1(M,g) + 2 \la \mu 
 \mu_1(M,g)  \int_M |v|^{N-2} v v_{\ep} \,dv_g}{ \la^2 \int_M |u_{\ep}|^{N-2}
v_{\ep}^2 \,dv_g  + \mu^2 \int_M |u_{\ep}|^{N-2} v^2 \,dv_g  + 2\la \mu 
\int_M |u_{\ep}|^{N-2} v v_{\ep} \,dv_g} \leq 1.$$
If $\la \mu < 0$ then, we write that since $N-2 \in ]0,1[$, 
$$ |u_{\ep}|^{N-2} \leq Y( v_{\ep}) v_{\ep}^{N-2} + \mu_1(M,g) v^{N-2}.$$
We obtain that 
\begin{eqnarray*}    
\la^2 \int_M |u_{\ep}|^{N-2}
v_{\ep}^2 \,dv_g  + \mu^2 \int_M |u_{\ep}|^{N-2} v^2 \,dv_g  + 2\la \mu 
\int_M |u_{\ep}|^{N-2} v v_{\ep} \,dv_g  \\
 \geq  \la^2 Y(v_{\ep})+ \mu^2 
 \mu_1(M,g)  - C \left(\int_M v^{N-1} v_{\ep} \,dv_g + \int_M v_{\ep}^{N-1} v
   \,dv_g \right).
\end{eqnarray*}
where $C >0$ is as in in the following a positive real number independent
of $\ep$. Together with (\ref{normp}), we get that 
$$\la^2 \int_M |u_{\ep}|^{N-2}
v_{\ep}^2 \,dv_g  + \mu^2 \int_M |u_{\ep}|^{N-2} v^2 \,dv_g  + 2\la \mu 
\int_M |u_{\ep}|^{N-2} v v_{\ep} \,dv_g  
 \geq  \la^2 Y(v_{\ep})+ \mu^2 \mu_1(M,g) + O(\ep^{\frac{n-2}{4}}).$$
It follows that 
\begin{eqnarray} \label{r1}
 \sup_{(\lambda,\mu) \in \mR^2 \setminus \{(0,0)\}}  \frac{  \la^2  Y(v_{\ep}) + \mu^2 \mu_1(M,g) + 2 \la \mu 
 \mu_1(M,g)  \int_M |v|^{N-2} v v_{\ep} \,dv_g}{ \la^2 \int_M |u_{\ep}|^{N-2}
v_{\ep}^2 \,dv_g  + \mu^2 \int_M |u_{\ep}|^{N-2} v^2 \,dv_g  + 2\la \mu 
\int_M |u_{\ep}|^{N-2} v v_{\ep} \,dv_g} \leq 1 +O(\ep^{\frac{n-2}{4}}).
\end{eqnarray}
By Lemma \ref{estim}, 
\begin{eqnarray*} 
 \int_M u_{\ep}^N  \,dv_g
& \leq &   
{(Y(u_{\ep}))}^{\frac{n}{2}} \int_M 
v_{\ep}^N \,dv_g + \mu_1(M,g)^{\frac{n}{2}} \int_M 
v^N \,dv_g \\ 
& & + C  \left(\int_M v^{N-1} v_{\ep} \,dv_g + \int_M v_{\ep}^{N-1} v
   \,dv_g \right).
\end{eqnarray*}
By (\ref{intv}), (\ref{intve}), (\ref{test}) and (\ref{normp}), we obtain
\begin{eqnarray} \label{r2}
{\left(  \int_M u_{\ep}^N  \,dv_g \right)}^{\frac{2}{n}}
\leq  {( \mu_1(M,g)^{\frac{n}{2}}+
  \mu_1(\mS^n)^{\frac{n}{2}})}^{\frac{2}{n}} -C \ep^2  + O(\ep^{\frac{n-2}{4}}) + o(\ep^2).
\end{eqnarray}
Since $\frac{n-2}{4} > 2$, we get from (\ref{r1}) and (\ref{r2}) that for
$\ep$ small enough
\begin{eqnarray*} 
\mu_2(M,g) &\leq  & \sup_{(\lambda,\mu) \in \mR^2 \setminus \{(0,0)\}}
F(u_{\ep}, \la v_{\ep} + \mu v) \\
&\leq & {( \mu_1(M,g)^{\frac{n}{2}}+
  \mu_1(\mS^n)^{\frac{n}{2}})}^{\frac{2}{n}} -C \ep^2  + O(\ep^{\frac{n-2}{4}}) +
o(\ep^2) <  {( \mu_1(M,g)^{\frac{n}{2}}+
  \mu_1(\mS^n)^{\frac{n}{2}})}^{\frac{2}{n}}.
\end{eqnarray*}
This proves Theorem \ref{upbound} if $\mu_1(M,g) >0$.

Now, we assume that $\mu_1(M,g)=0$, that $(M,g)$ is not locally conformally
flat  and that $n \geq 9$. 
For more simplicity, We set $u_{\ep} =
v_{\ep}$  instead of  $ u_{\ep} = Y(v_{\ep})^{\frac{n-2}{4}}  v_{\ep}$ as
above. We proceed exactly  as in the case $\mu_1(M,g)>0$. We obtain that 
for $(\lambda,\mu) \in \mR^2 \setminus \{(0,0)\}$ 
\begin{eqnarray*}
F(u_{\ep}, \la v_{\ep} + \mu v) & = &
 \frac{  \la^2  Y(v_{\ep})}{ \la^2 \int_M v_{\ep}^N
\,dv_g  + \mu^2 \int_M |v_{\ep}|^{N-2} v^2 \,dv_g  + 2\la \mu 
\int_M |v_{\ep}|^{N-1} v \,dv_g} {\left( \int_M v_{\ep}^N  \,dv_g
\right)}^{\frac{2}{n}} \\
& =&  \frac{  \la^2  Y(v_{\ep}) }{ \la^2   + \mu^2 \int_M
 |v_{\ep}|^{N-2} v^2 \,dv_g   + 2\la \mu 
\int_M |v_{\ep}|^{N-1} v \,dv_g}.
\end{eqnarray*}

 Let $\la_{\ep},
\mu_{\ep}$ be such that $\la_{\ep}^2+ \mu_{\ep}^2 = 1$ and such that
$$ F(u_{\ep}, \la_{\ep}  v_{\ep} + \mu_{\ep} v) =  \sup_{(\lambda,\mu) \in
  \mR^2 \setminus \{(0,0)\}} (u_{\ep},  \la v_{\ep} + \mu  v).$$
 If $\la_{\ep} = 0$,
we obtain that $F(u_{\ep}, \la_{\ep} v_{\ep} + \mu_{\ep} v) = 0$ and the
theorem would be proven. 
Then we assume that $\la_{\ep} \not= 0 $ and we write that 
$$F(u_{\ep}, \la_{\ep} v_{\ep} + \mu_{\ep} v) = 
\frac{Y(v_{\ep}) }{ 1 + 2 x_{\ep} b_{\ep} + x_{\ep}^2 a_{\ep}}$$
where $x_{\ep} = \frac{\mu_{\ep}}{\la_{\ep}}$ and where, using (\ref{normp})
$$b_{\ep} = \int_M v_{\ep}^{N-1} v \,dv_g \sim_{\ep \to 0} C \ep^{\frac{n-2}{4}}$$
and 
$$a_{\ep} =  \int_M v_{\ep}^{N-2} v^2 \,dv_g \sim_{\ep \to 0} C \ep.$$
Maximizing this expression in $x_{\ep}$ and using (\ref{test}), we get that 
$$F(u_{\ep}, \la_{\ep} v_{\ep} + \mu_{\ep} v) \leq \frac{\mu_1(\mS^n) - C(M)
  \ep^2+ o(\ep^2)}{1 - \frac{b_{\ep}^2 }{a_{\ep}}} =  \frac{\mu_1(\mS^n) - C(M)
  \ep^2+ o(\ep^2)}{1 - O(\ep^{\frac{n-4}{2}})}.$$
Since $n \geq 9$, $\frac{n-4}{2} > 2$ and we get that for $\ep$ small,
$$F(u_{\ep}, \la_{\ep} v_{\ep} + \mu_{\ep} v)< \mu_1(\mS^n).$$
This proves Theorem \ref{upbound}.

\section{Existence of a minimum of $\mu_2(M,g)$} \label{sectatt}

The aim of this section is to prove Theorem~\ref{attain}.

\setcounter{step}{0}
We study a sequence of metrics 
$(g_m)_m= (u_m ^{N-2} g)_m$ 
($u_m >0$, $u_m \in C^{\infty}(M)$) which minimizes the infimum
in the definition of $\mu_2(M,g)$ i.e. a sequence of metrics such that 
$$\lim_m \lambda_2(g_m) {\Vol(M,g_m)}^{\frac{2}{n}}=\mu_2(M,g).$$
Without loss of generality, we may assume that $\Vol(M,g_m)=  1$ i.e. that 
\begin{eqnarray} \label{u_bound} 
 \int_M u_m^N \,dv_g =1.
\end{eqnarray}
In particular, the sequence $(u_m)_m$ is bounded in $L^N(M)$ and there
exists $u \in L^N(M)$, $u \geq 0$ 
such that $u_m \weakto u$ weakly in $L^N(M)$. We are
going to prove that $u \not=0 $ and that the generalized metric $u^{N-2} g$
minimizes $\mu_2(M,g)$.
Proposition \ref{la1la2} implies the existence of  
$v_m, w_m \in C^{\infty}(M)$, $v_m \geq 0$ 
such that 
\begin{eqnarray} \label{eqvm}
L_g v_m = \la_{1,m}   u_m^{N-2} v_m
\end{eqnarray}
and
\begin{eqnarray} \label{eqwm}
L_g w_m = \la_{2,m}    u_m^{N-2} w_m.
\end{eqnarray}
where $\la_{i,m} =  \la_i (g_m)$ and such that  
\begin{eqnarray} \label{vw_bound}
\int_M u_m^{N-2} v_m^2 \,dv_g = \int_M u_m^{N-2} w_m^2 \,dv_g
= 1 \; \hbox{and} \; \int_M u_m^{N-2} v_m w_m  \,dv_g = 0
\end{eqnarray}
With these notations and by
(\ref{u_bound}),
 
$$\lim_{m}  \la_{2,m} = \mu_2(M,g).$$

 Moreover, by the maximum principle, $v_m >0$. If $\la_{1,m} =
\la_{2,m}$ then $w_m$ would be a minimizer of the functional associated to
$\la_{1,m}$ and by the maximum principle, we would get that 
$w_m >0$. This contradicts (\ref{vw_bound}). Hence,   
$\la_{1,m} <  \la_{2,m}$ for all $m$.
 The sequences $(v_m)_m$ and $(w_m)_m$ are bounded in
$H_1^2(M)$. We can find $v,w \in H_1^2(M)$, $v \geq 0$ such that 
$v_m$ (resp. $w_m$) tends to $v$ (resp. $w$) weakly in $
H_1^2(M)$. Together with the weak convergence of the $(u_m)_m$ towards $u$
in $L^N(M)$, we get that in the sense of distributions

\begin{eqnarray} \label{eqvlim}
L_g v  =  \widehat\mu_1  u^{N-2} v
\end{eqnarray}
and
\begin{eqnarray} \label{eqwlim}
L_g w = \mu_2(M,g)\, u^{N-2} w.
\end{eqnarray}
where $\widehat\mu_1 = \lim_m \la_{1,m} \leq \mu_2(M,g)$. 

{}From what we know until now, it is not clear whether $v$ and $w$ are 
linearly independent, and even if they are, their restrictions to the set
$M\setminus u^{-1}(0)$ might be linearly dependent.
 
\noindent It will take a certain effort to prove the following claim. 

\begin{claim}
The functions 
$u^{N-2\over 2} v$ and  $u^{N-2\over 2} w$ are linearly independent.
\end{claim}

\noindent Once the claim is proved, we have $\spann(v,w)\in \Gr2u{H_1^2(M)}$, and this
implies that 
  $$\sup_{(\la,\mu)\neq (0,0)}  F(u, \la v + \mu w)  = \mu_2(M,g).$$ 
Hence, by equations (\ref{eqvlim}) and (\ref{eqwlim}), the generalized metric $u^{N-2} g$ minimizes $\mu_2(M,g)$, i.e. 
Theorem~\ref{attain} is proved.

The first step in the proof of the claim is an estimate that avoids 
concentration of $w_m$ and $v_m$.

\begin{step} Let $x \in M$ and $\ep \in \mathopen]0, \frac{N-2}{2}\mathclose[$. 
We choose a  cut-off function $\eta \in C^{\infty}$ such that $0 \leq \eta \leq 1$, $\eta(B_x(\de)) \equiv 1$ (where
$\delta>0$ is a small number) and $\eta(M \setminus B_x(2\de) ) \equiv 0$, $|\na \eta|\leq 2/\delta$. 
We define $W_m =  \eta |w_m|^{\ep} w_m$. Then, we have 
\begin{eqnarray} \label{inegstep1}
{\left( \int_M |W_m|^N
    \,dv_g \right)}^{\frac{2}{N}} \leq \mu_2(M,g) (1 - \al_{\ep})^{-1} \mu_1(\mS^n)^{-1}
{\left( \int_{B_x(2 \de)} u_m^N \right)}^{\frac{2}{n}}  {\left( \int_M |W_m|^N
    \,dv_g \right)}^{\frac{2}{N}} +C_{\de}.
\end{eqnarray}
where $C_{\de}$ is a constant that may depend on $\de$ but not on
$\ep$ and where  $\lim_{\ep \to 0} \al_{\ep}=0$.
Moreover, the same conclusion is true with $V_m = \eta |v_m|^{\ep}
v_m$ instead of $W_m$.
\end{step}
The proof uses classical methods. 
We will explain the proof for $W_m$. The proof
for $V_m$ uses exactly the same arguments.

At first, we differentiate the definition of $W$ and obtain
\begin{eqnarray} 
|\nabla W_m|^2 & \geq & \Bigl|\nabla( |w_m|^{\ep} w_m )\Bigr|^2 \eta^2 
- \bigl(2 |\nabla \eta|\, |w_m|^{1+\ep}\bigr)\left( \Bigl|\nabla  (|w_m|^{\ep} w_m) \Bigr| 
\eta\right) +  |\nabla \eta|^2 \, |w_m|^{2+2\ep} \nonumber\\
& \geq &  \Bigl|\nabla \bigl(|w_m|^{\ep} w_m\bigr) \Bigr|^2  \eta^2 
- \left( \frac{1}{2}  \Bigl|\nabla
 ( |w_m|^{\ep} w_m) |^2  \eta^2   + 2  |\nabla
\eta|^2 |w_m|^{2+2\ep} \right) +  |\nabla
\eta|^2 |w_m|^{2+2\ep}\nonumber
\end{eqnarray}

This leads to
\begin{equation}
 \eta^2  |\nabla (|w_m|^{\ep} w_m )|^2  \leq 2 |\nabla W_m|^2 + 2 
 |\nabla
\eta|^2 |w_m|^{2+2\ep}.
\label{ineq.etana}
\end{equation}

Now, we  want to derive lower bound for  
\begin{eqnarray} \label{W1}
( \nabla (\eta^2 |w_m|^{2 \ep} w_m), \nabla w_m )  = 
 | \nabla W_m|^2 - \bigl|\nabla (\eta |w_m|^{\ep })\bigr|^2\, |w_m|^2
\end{eqnarray}

For the second summand  on the right hand side in \eref{W1} we have the bound

\begin{eqnarray*} 
|\nabla (\eta |w_m|^{\ep})|^2 |w_m|^2  &  =  &|\nabla \eta|^2 |w_m|^{2+2
  \ep}+ 2 (\nabla \eta, \nabla |w_m|^{\ep}) \,\eta |w_m|^{2+\ep} +  
   \eta^2 \Bigl|\nabla (|w_m|^{\ep})\Bigr|^2 w_m^2
\nonumber \\
& \leq & 2 |\nabla \eta|^2 |w_m|^{2+2
  \ep} + 2  \eta^2\Bigl|\nabla (|w_m|^{\ep})\Bigr|^2 w_m^2\\ 
& \leq &  2 |\nabla \eta|^2 |w_m|^{2+   2 \ep} 
+\frac{2  \eta^2\ep^2}{(1+\ep)^2} \Bigl|\nabla
(|w_m|^{\ep} w_m) \Bigr|^2\\
& \leq & (2 +\frac{4 \ep^2}{(1+\ep)^2}) 
 |\nabla \eta|^2 |w_m|^{2+   2 \ep}+ \frac{4 \ep^2}{(1+\ep)^2} |\nabla
 W_m|^2.
\end{eqnarray*}
Here, we used \eref{ineq.etana} in the last line.
Coming back to (\ref{W1}), we obtain that 
$$ ( \nabla (\eta^2 |w_m|^{2 \ep} w_m), \nabla w_m )  \geq (1-\al_{\ep}) 
 |\nabla
 W_m|^2 - C  |\nabla \eta|^2 |w_m|^{2+   2 \ep}. $$
where $\al_{\ep} \to 0$ when $\ep \to 0$ and where $C>0$ is a constant
independent of $\ep$.
This relations shows that

\begin{eqnarray*} 
\int_M \eta^2 |w_m|^{2\ep} w_m  L_g(w_m)  \,dv_g  \geq (1-\al_{\ep})\int_M
c_n |\nabla W_m|^2
\,dv_g -C \int_M  |\nabla \eta|^2 |w_m|^{2+   2 \ep} \,dv_g+ 
\min\Scal \int W_m^2\, dv_g.
\end{eqnarray*}
Now, since $\ep < \frac{N-2}{2}$,the sequence $(w_m)_m$ is bounded in
$L^{2+2\ep}(M)$ (and hence the sequence $(W_m)_m$ is bounded in
$L^2(M)$). As a consequence, there exists a constant
$C_{\de}$ possibly depending on~$\de$ but not on~$\ep$, and such that 
\begin{eqnarray} \label{W4} 
\int_M \eta^2 |w_m|^{2\ep} w_m  L_g(w_m)  \,dv_g  \geq (1-\al_{\ep})\int_M
\left(c_n |\nabla W_m|^2 + B_0(M,g) W_m^2\right) \,dv_g - C_{\de}.
\end{eqnarray}
Using  equation (\ref{eqwm}) in the left hand side of (\ref{W4}) and
applying Sobolev inequality $(S)$ to the right hand side, we get that 
\begin{eqnarray*} 
\mu_2(M,g) 
\int_M u_m^{N-2} W_m^2 \,dv_g \geq (1 - \al_{\ep}) \mu_1(\mS^n) {\left( \int_M |W_m|^N
    \,dv_g \right)}^{\frac{2}{N}}  - C_{\de}.
\end{eqnarray*}
By the H\"older inequality, we obtain
\begin{eqnarray*}
{\left( \int_M |W_m|^N
    \,dv_g \right)}^{\frac{2}{N}} \leq \mu_2(M,g) (1 - \al_{\ep})^{-1} \mu_1(\mS^n)^{-1}
{\left( \int_{B_x(2 \de)} u_m^N \right)}^{\frac{2}{n}}  {\left( \int_M |W_m|^N
    \,dv_g \right)}^{\frac{2}{N}} +C_{\de}.
\end{eqnarray*}
This ends the proof of the step.
\begin{step} 
 If $\mu_2(M,g) < \mu_1(\mS^n)$, then the generalized metric $u^{N-2} g$
 minimizes $\mu_2(M,g)$.
\end{step}
{}From (\ref{inegstep1}), and the fact  $\mu_2(M,g) < \mu_1(\mS^n)$, we get
that for $\ep$ small enough, there exists a constant $K <1$ such that

$${\left( \int_M |W_m|^N
    \,dv_g \right)}^{\frac{2}{N}} \leq  K 
{\left( \int_{B_x(2 \de)} u_m^N \right)}^{\frac{2}{n}}  {\left( \int_M |W_m|^N
    \,dv_g \right)}^{\frac{2}{N}} +C_{\de}.$$
Since $\int_{B_x(2 \de)} u_m^N   \leq 1$,
the sequence $\int_M |W_m|^N \,dv_g $ is bounded.  
This implies that $(w_m)_m$ is bounded in 
$L^{N+{\ep}}(B_x(\de))$ and since $x$ is arbitrary in
$L^{N+\ep}(M)$. Weak convergences $w_m\to w$ in $H_1^2(M)$ implies
strong convergence $w_m\to w$ in $L^{N-\ep}(M)$. The H\"older inequality
yields then strong convergence in $L^N(M)$. 
After passing to a subsequence we obtain that  $(w_m)_m$ tends to $w$
strongly in $L^N (M)$. This implies that we can pass to the limit in
(\ref{vw_bound}) and hence that 
$u^{\frac{N-2}{2}} v$ and  $u^{\frac{N-2}{2}} w$ are
linearly independent. The claim follows in this case.\\

In the following, we assume that $\mu_1(M,g) >0$ and that 
$$\mu_2(M,g) < {\left( \mu_1(M,g)^{\frac{n}{2}} + \mu_1(\mS^n)^{\frac{n}{2}} 
\right)}^{\frac{2}{n}}.$$
We define the \emph{set of concentration points} 
$$\Om = \Biggl\{ x \in M \,\Big|\, \forall \de>0, \;  \limsup_m  \int_{B_x(\de)}
  u_m^N \,dv_g > \frac{1}{2} \Biggr\}.$$ 
Since $\int_M u_m^N \,dv_g = 1$, we can assume --- after passing to
a subsequence --- that $\Om$ contains at most one point.

We now prove that:

\begin{step} \label{st2}
Let $U$ be an open set such that $\overline{U}  \subset M \setminus
\Om$. Then, the sequence $(v_m)_m$ (and $(w_m)_m$ resp.) converges 
towards $v$ (and $w$ resp.) strongly in $H_1^2(\overline{U})$.
\end{step}
Without loss of generality, we prove the result only for $w$. 
For any $x \in M \setminus \Om$ we can find $\de>0$
with
  $$\limsup_m \int_{B_x(2 \de)} u_m^N \,dv_g\leq{1\over 2}.$$
Using 
$\mu_2(M,g) < {(\mu_1(M,g)^{\frac{n}{2}} +
\mu_1(\mS^n)^{\frac{n}{2}})}^{\frac{2}{n}} \leq 2^{\frac{2}{n}} \mu_1(\mS^n)$
we obtain
for a small $\ep>0$
$$ \mu_2(M,g) (1 - \al_{\ep})^{-1} \mu_1(\mS^n)^{-1}
{\left( \int_{B_x(2 \de)} u_m^N \right)}^{\frac{2}{n}} \leq K <1$$
for almost all $m$.
Together with inequality (\ref{inegstep1}), this proves that $\int_M
|W_m|^N \,dv_g $ is bounded. This implies that $(w_m)_m$
is bounded in $L^{N+{\ep}}(B_x(\de))$. As in last step, this proves that up
to a subsequence,
$(w_m)_m$ tends to $w$ strongly in $L^N(U)$. Using equation (\ref{eqwm})
and (\ref{eqwlim}), we easily obtain that 
$$\lim_m \int_U |\nabla w_m|^2 dv_g =\int_U |\nabla w|^2 dv_g.$$  
Together with the weak convergence of $(w_m)_m$ to $w$, this proves the step.\\

 Now, we set for all $m$,
$$S_m= \{ \lambda v_m + \mu w_m | \la^2+\mu^2 = 1 \}  \; \hbox{ and } \;
S=  \{ \lambda v + \mu w | \la^2+\mu^2 = 1 \}.$$

\begin{step} \label{st1}
There exists a sequence $(\w_m)_m$ ($\w_m \in S_m$)  and
$\w \in S$ such that $\w_m$ tends to $\w$
strongly in $H_1^2(M)$.
\end{step}

By theorem \ref{sobolev}, there exists $\la_m,\mu_m$ such that 
$\la_m^2+\mu_m^2= 1$ and 
 such that 
\begin{align} 
2^{2/n} \mu_1(\mS^n) &\int_M
u_m^{N-2} {(\la_m(v_m-v) + \mu_m(w_m - w))}^2 \,dv_g\nonumber  \\
 \leq      &\int_M c_n |\nabla(\la_m(v_m-v)
 +   \mu_m(w_m - w))   |^2 \,dv_g \label{sob}  \\ 
{}+  & \int_M  B_0(M,g)(\la_m(v_m-v) + \mu_m(w_m - w))^2 \,dv_g  \nonumber 
\end{align} 
Up to a subsequence, there exists $\la,\mu$ such that $\la^2+\mu^2=1$ and
such that $\lim_m \la_m = \la$ and $\lim_m \mu_m = \mu$. We set  
$\w_m= \la_m v_m + \mu_m w_m \in  S_m$ and  $\w= \la v
+ \mu w$.  Then, $\w_m$ tends to $\w$ weakly in $H_1^2(M)$.
A first remark is that by strong convergence in $L^2(M)$   
\begin{eqnarray}  \label{sob1}
\lim_m \int_M (\la_m(v_m-v) + \mu_m(w_m - w))^2 \,dv_g  = 0.
\end{eqnarray}
Using the weak convergence of $\w_n$ to $\w$ in  $H_1^2(M)$ and the weak
convergence of $u_m$ to $u$ in $L^N(M)$,  
it is easy to compute that 
\begin{eqnarray} \label{sob2}
\int_M
u_m^{N-2} {(\la_m(v_m-v) + \mu_m(w_m - w))}^2 \,dv_g   = \int_M u_m^{N-2}
\w_m^2 \,dv_g - \int_M u^{N-2}
\w^2 \,dv_g + o(1)
\end{eqnarray}
and that 
\begin{eqnarray*}
\int_M  c_n |\nabla(\la_m(v_m-v) + \mu_m(w_m - w))   |^2 \,dv_g 
&  =  &
\la^2 \left( \int_M c_n |\nabla v_m |^2 \,dv_g - \int_M c_n |\nabla v |^2 \,dv_g
  \right)\\
& + &  \mu^2  \left( \int_M c_n |\nabla w_m |^2 \,dv_g - \int_M
    c_n |\nabla w |^2 \,dv_g \right)\\
&+ &   
2 \la \mu \left( \int_M c_n (\nabla v_m, \nabla w_m)  \,dv_g - \int_M  c_n
  (\nabla v, \nabla w)  \,dv_g \right)  +o(1).
\end{eqnarray*} 
Using equations (\ref{eqvm}), (\ref{eqwm}), (\ref{eqvlim}) and
(\ref{eqwlim}),
we get that 
\begin{eqnarray*}
 \int_M  c_n |\nabla(\la_m(v_m-v) + \mu_m(w_m - w))   |^2 \,dv_g 
&  =  & 
 \la^2\,\widehat\mu_1 \left(  \int_M u_m^{N-2} v_m^2 \,dv_g - 
\int_M u^{N-2} v^2) \,dv_g \right)\\
&  + &  \mu^2\, \mu_2(M,g)  \left(  \int_M u_m^{N-2} w_m^2 \,dv_g - 
\int_M u^{N-2} w^2) \,dv_g \right) \\
& + & 2 \la \mu\, \mu_2(M,g)  \left(  \int_M u_m^{N-2} v_ m w_m \,dv_g -
 \int_M u^{N-2} v w \,dv_g \right) + o(1).
\end{eqnarray*} 
Since $\widehat\mu_1 \leq \mu_2(M,g)$ and since, by weak convergence  
$$\liminf_m  \int_M u_m^{N-2} v_m^2 \,dv_g - 
\int_M u^{N-2} v^2) \,dv_g  \geq 0,$$
we get that 
\begin{eqnarray*}
\int_M c_n |\nabla(\la_m(v_m-v) + \mu_m(w_m - w))   |^2 \,dv_g 
& \leq  &  
    \la^2\,\mu_2(M,g) \left(  \int_M u_m^{N-2} v_m^2 \,dv_g - 
    \int_M u^{N-2} v^2 \,dv_g \right)\\
& +  & 
    \mu^2\,\mu_2(M,g) \left(  \int_M u_m^{N-2} w_m^2 \,dv_g - 
    \int_M u^{N-2} w^2 \,dv_g \right)  \\
& +&  
    2 \la \mu\, \mu_2(M,g)  \left(  \int_M u_m^{N-2} v_ m w_m \,dv_g -
    \int_M u^{N-2} v w \,dv_g \right), 
\end{eqnarray*}
and hence, 
\begin{eqnarray} \label{sob3} 
\int_M c_n |\nabla(\la_m(v_m-v) + \mu_m(w_m - w))   |^2 \,dv_g  \leq  \mu_2(M,g)  \left( \int_M u_m^{N-2} \w_m^2 \,dv_g - 
    \int_M u^{N-2} \w^2 \,dv_g \right) +o(1). 
\end{eqnarray}
Together with (\ref{sob}), (\ref{sob1}) and (\ref{sob2}), we obtain that 
\begin{eqnarray*}
&&2^{2/n} \mu_1(\mS^n) \left( \int_M u_m^{N-2} \w_m^2 \,dv_g - 
\int_M u^{N-2} \w^2 \,dv_g \right)  \\
&&\leq  \mu_2(M,g)  \left( \int_M u_m^{N-2} \w_m^2 \,dv_g - 
\int_M u^{N-2} \w^2 \,dv_g \right) +o(1).
\end{eqnarray*}

 Since $\mu_2(M,g) < {(\mu_1(M,g)^{\frac{n}{2}} +
  \mu_1(\mS^n)^{\frac{n}{2}})}^{\frac{2}{n}} \leq 2^{\frac{2}{n}} \mu_1(\mS^n)$, we get
that 
$$\left( \int_M u_m^{N-2} \w_m^2 \,dv_g - 
\int_M u^{N-2} \w^2 \,dv_g \right)  \leq K_0   \left( \int_M u_m^{N-2} \w_m^2 \,dv_g - 
\int_M u^{N-2} \w^2 \,dv_g \right) +o(1)$$
where $K_0 <1$. This implies that
\begin{eqnarray}  \label{unot=0}
1= \lim_m \int_M u_m^{N-2} \w_m^2 \,dv_g = \int_M u^{N-2} \w^2 \,dv_g
\end{eqnarray}
and hence by (\ref{sob3}). 
$$\lim_m \int_M c_n |\nabla(\la_m(v_m-v) + \mu_m(w_m - w))   |^2 \,dv_g =0.$$
The step easily follows.\\

As a remark, (\ref{unot=0}) implies that  $u^{\frac{N-2}{2}} \w
\not\equiv 0$.

 Now, we set $\v_m = - \mu_m v_m + \la_m w_m$ and  $\v= -\mu v +
\la w$. We prove that 

\begin{step} \label{st3}
There exists $x \in M$ such that  
$$\limsup_m \int_{B_x{\delta}}  u_m^2 (\v_m- \v)^2 \,dv_g = 1$$
for all $\delta>0$.
\end{step}
 
The sequence 
$(\v_m)_m$ tends to $\v$ weakly in $H_1^2(M)$. 
If $\Om=\emptyset$, then we know from Step \ref{st2} that 
$(\v_m)_m$ tends to $\v$ strongly in $H_1^2(M)$, which implies
$\int u^{N-2}\bar v \bar w=0$. Hence, in the case $\Omega=\emptyset$,
the functions  $ u^{N-2\over 2}\bar v$ and $ u^{N-2\over 2}\bar w$  
are linearly independent, and the claim follows.

Hence, without loss of generality let $\Om = \{ x \}$ where $x$ 
is some point of $M$. We assume that the claim is false, i.e.  
$u^{\frac{N-2}{2}} v$ and $ u^{\frac{N-2}{2}}w$ are 
linearly dependent. As $u^{\frac{N-2}{2}} \w\not\equiv 0$,
there exists  $b\in \mR$ with $u^{\frac{N-2}{2}} \v= b u^{\frac{N-2}{2}} \w$.
Hence, 
$$ 0 = \int_M u^{N-2} \v^2 \,dv_g + b^2\int_M u^{N-2} \w^2 \,dv_g - 
2 b \int_M u^{N-2} \v\, \w \,dv_g.$$
By strong convergence of $(\w_m)_m$ to $\w$ in $H_1^2(M)$, weak 
convergence of $(\v_m)_m$ to $\v$ in $H_1^2(M)$ and weak convergence of $(u_m)_m$ to $u$ in
 $L^N(M)$, we have $ \int_M u^{N-2} \w^2 \,dv_g= 1 $ and 
 $\int_M u^{N-2} \v \, \w \,dv_g = 0$.
We obtain $\int_M u^{N-2} \v^2 \,dv_g + b^2 = 0$. As a consequence,
$u^{\frac{N-2}{2}} \v \equiv 0$. Let now $\delta>0$. We write that
\begin{eqnarray*}
\int_{B_{x}(\de)}   u_m^{N-2}( \v_m-\v)^2 \,dv_g
&  = &  \int_{B_{x}(\de)}   u_m^{N-2}\v_m^2 \,dv_g\\
& = & 1 - \int_{M \setminus B_{x}(\de)}   u_m^{N-2} \v_m^2 \,dv_g.
\end{eqnarray*}
  By step \ref{st2}, 
$$\lim_m \int_{M \setminus B_{x}(\de)}   u_m^{N-2} \v_m^2 \,dv_g= \int_{M
  \setminus B_{x}(\de)}  u^{N-2} \v^2 \,dv_g =0.$$
This proves the step.

\begin{step} \label{st4}
Conclusion.
\end{step}

Let $\de >0$ be a small fixed number. In the following, $o(1)$
denotes a sequence of real numbers which tends to $0$, however we do not
claim that the convergence is uniform in $\delta$. 
By step \ref{st3} and the H\"older inequality,
\begin{eqnarray*}
1 & = & \int_{B_x(\de)}   u_m^{N-2} (\v_m- \v)^2 \,dv_g + o(1)\\ 
  & \leq & {\left( \int_{B_x(\de)}  u_m^N \,dv_g
  \right)}^{\frac{2}{n}} {\left( \int_M |\v_m -\v|^N 
  \,dv_g  \right)}^{\frac{2}{n}}+o(1).
\end{eqnarray*}
Applying Sobolev inequality $(S)$, we get that 
\begin{eqnarray*}
1 & \leq  & {\left( \int_{B_x(\de)}  u_m^N \,dv_g
  \right)}^{\frac{2}{n}}  \mu_1(\mS^n)^{-1} \left( \int_M c_n 
|\nabla (\v_m - \v) |^2 \,dv_g +B_0(M,g) \int_M (\v_m - \v)^2 \,dv_g \right) + o(1).
\end{eqnarray*}
By strong convergence of $(\v_m-\v)_m$ to $0$ in $L^2(M)$, 
\begin{eqnarray*}
1& \leq  &   {\left( \int_{ B_x(\de)}  u_m^N \,dv_g
  \right)}^{\frac{2}{n}}  \mu_1(\mS^n)^{-1} \left( \int_M c_n 
|\nabla (\v_m - \v) |^2+  S_g  (\v_m - \v)^2 \,dv_g \right)+ o(1) 
\end{eqnarray*}
Using equations   (\ref{eqvm}), (\ref{eqwm}), (\ref{eqvlim}),
(\ref{eqwlim}) and the fact that $\widehat\mu_1 \leq \mu_2(M,g)$,
we get that 
\begin{eqnarray*}
1& \leq &  {\left( \int_{ B_x(\de)}  u_m^N \,dv_g
  \right)}^{\frac{2}{n}} \mu_1(\mS^n)^{-1} \mu_2(M,g) 
   \int_M  u_m^{N-2}( \v_m -  \v)^2 \,dv_g  \\
& =& {\left( \int_{B_x(\de)}  u_m^N \,dv_g
  \right)}^{\frac{2}{n}} \mu_1(\mS^n)^{-1} \mu_2(M,g).
\end{eqnarray*}
Since $\mu_2(M,g) < {(\mu_1(M,g)^{\frac{n}{2}} +
  \mu_1(\mS^n)^{\frac{n}{2}})}^{\frac{2}{n}}$, we obtain that  
\begin{eqnarray*}
\int_{ B_x(\de)}  u_m^N \,dv_g > \frac{\mu_1(\mS^n)^{\frac{n}{2}}}
{\mu_1(M,g)^{\frac{n}{2}} +
  \mu_1(\mS^n)^{\frac{n}{2}}}.
\end{eqnarray*}
and since $\int_M u_m^N \,dv_g = 1$, 
\begin{eqnarray} \label{in_final}
\int_{M \setminus B_x(\de)}  u_m^N \,dv_g < \frac{\mu_1(M,g)^{\frac{n}{2}}}
{\mu_1(M,g)^{\frac{n}{2}} +
  \mu_1(\mS^n)^{\frac{n}{2}}}.
\end{eqnarray}
Now, we write that by strong convergence of $(\w_m)_m$ in $H_1^2(M)$,
\begin{eqnarray*}
a_{\de} & = & 
\int_{B_x(\de)}   u_m^{N-2} \w_m^2 \,dv_g 
\end{eqnarray*}

\begin{eqnarray*}
1 -a_{\de} & = & 
\int_{M \setminus  B_x(\de)}   u_m^{N-2} \w_m^2 \,dv_g 
\end{eqnarray*}
where $a_{\de}$ does not depend of $m$ and tends to $0$ when $\de$ tends to
$0$. By H\"older inequality,
\begin{eqnarray*}
1 -a_{\de} & \leq & {\left( \int_{M \setminus B_x(\de)}  u_m^N \,dv_g
  \right)}^{\frac{2}{n}} {\left( \int_M \w^N 
\,dv_g  \right)}^{\frac{2}{n}}.
\end{eqnarray*}
Since $\mu_1(M,g)$ is the minimum of Yamabe functional, we get that 
\begin{eqnarray*}
1 -a_{\de} & \leq & {\left( \int_{M \setminus B_x(\de)}  u_m^N \,dv_g
  \right)}^{\frac{2}{n}} \mu_1(M,g)^{-1} \int_M  \left(c_n 
|\nabla \w_m|^2 + S_g  \w_m^2\right) \,dv_g  .
\end{eqnarray*}
 As we did for $\v$, we obtain

\begin{eqnarray*}
1 -a_{\de} & \leq & {\left( \int_{M \setminus B_x(\de)}  u_m^N \,dv_g
  \right)}^{\frac{2}{n}} \mu_1(M,g)^{-1} \mu_2(M,g)
  \underbrace{\int_M    u_m^{N-2} \w_m^2 \,dv_g}_1 
\end{eqnarray*}
By (\ref{in_final}), in the limit $\de \to 0$, this gives 
$$\mu_2(M,g) \geq {(\mu_1(M,g)^{\frac{n}{2}} +
  \mu_1(\mS^n)^{\frac{n}{2}})}^{\frac{2}{n}}.$$
This is false by assumption. Hence, the claim is proved, and 
Theorem~\ref{attain} follows.

\section{The invariant $\mu_k(M)$ for $k \geq 3$}
A natural question is: Can we do the same work for $\mu_k(M)$ with $k \geq
3$? This problem is still open but seems to be hard. Let $(M,g)$ be a
compact Riemannian manifold of dimension $n \geq 3$.  
Using the variational
characterization of $\mu_k(M)$, one can check that $\mu_k(M) \leq
k^{\frac{2}{n}} \mu_1(\mS^n)$. It is natural to conjecture that one has
equality if $M$ is the round sphere i.e. that 
$\mu_k(\mS^n)= k^{\frac{2}{n}} \mu_1(\mS^n)$. 
However, the following
result shows that is false:

\begin{prop} \label{muk}
Let $n \in \mN^*$. Then, for $n\geq 7$ 
  $$\mu_{n+2}(\mS^n) < (n+2)^{\frac{2}{n}} \mu_1(\mS^n).$$
\end{prop}

{\bf Proof:} Let us study $\mS^n$ with its natural embedding into
$\mR^{n+1}$. We have $L_g(1)=n(n-1)$. Hence, $\la_1(\mS^n) \leq n(n-1)$. Let
also $x_i$ ($i \in [1,\cdots, n+1]$) be the  canonical
coordinates on $\mR^{n+1}$. As one can check,
$$L_g(x_i) = \frac{n(n-1)(n+2)}{n-2} x_i$$
and hence 
$\la_{n+2}(\mS^n) \leq \frac{n(n-1)(n+2)}{n-2}$.
This shows that 
$$\mu_{n+2}(\mS^n) \leq \frac{n(n-1)(n+2)}{n-2} \om_{n}^{\frac{2}{n}}.$$
As one can check, for $n\geq 7$
$$\frac{n(n-1)(n+2)}{n-2} \om_{n}^{\frac{2}{n}} < (n+2)^{\frac{2}{n}}
n (n-1) \om_{n}^{\frac{2}{n}}=  (n+2)^{\frac{2}{n}}  \mu_1(\mS^n).$$
This ends the proof of Proposition \ref{muk}.

\section{The case of manifolds whose Yamabe invariant is negative} 
\label{negative}
We let $(M,g)$ be a compact Riemannian manifold of dimension $n \geq 3$. Then, we
have: 

\begin{prop}
Let $k \in \mN^*$. Assume that $\mu_k(M,g) < 0$. Then, $\mu_k(M,g)= - \infty$.
\end{prop}
\noindent {\bf Proof:} After a possible change of metric in the conformal
class, we can assume that $\la_k(g) <0$. This implies that we can find some
smooth functions $v_1, \cdots,v_k$ satisfying 
$$L_g v_i = \la_i(g) v_i $$
for all $i \in \{1, \cdots,k \}$ and such that 
$$\int_M v_i v_j dv_g = 0$$
 for all $i,j \in \{1, \cdots, k\}$, $i \not = j$. Let $v_\ep$ be defined
 as in the proof of Theorem \ref{upbound}. We define $u_{\ep} = v_{\ep} +
 \ep$ to obtain a positive function. We set 
$V= \{v_1, \cdots, v_k \}$.  It is easy to check that, uniformly in $v \in V$
$$\lim_{\ep to 0} \int_M v_{\ep}^{N-2} v^2 dv_g =0.$$ 
Since $\la_i <0$, it is then easy to 
see that $\sup_{v \in V} F(v_{\ep},v) = -\infty$. Together with the
variational characterization of $\mu_k(M,g)$, we get that $\mu_k(M,g) = -
\infty$.\\

\noindent This result proves for example that if the Yamabe invariant of
$(M,g)$ is negative, then $\mu_1(M,g) = - \infty$.
This is the reason why we restricted in this article to the case
of non-negative Yamabe invariant. Many of our results and proofs
remain valid in the case $\mu_2(M) \geq 0$. 
However, if the Yamabe invariant of $(M,g)$ is non-positive, there
are other ways  to find
nodal solutions of Yamabe equation. Indeed, Aubin's methods \cite{aubin:76}
 can be applied
to avoid concentration phenomenom. See for example    
\cite{djadli.jourdain:02},
\cite{jourdain:99}, \cite{holcman:99} for such methods.
Here, we present very briefly one new method in this case. We just sketch
it since it is not the purpose of our paper to find solutions of
Yamabe equation with Aubin's type methods.
 
At first, for any metric $\tilde{g}$ conformal to $g$,  
we let $\la_1^+(\tilde{g})$ be the first \emph{positive} eigenvalue of Yamabe
operator. We then define $\la^+ = \inf \la_1^+(\tilde{g})
\Vol(M,\tilde{g})^{\frac{2}{n}}$ where the infimum is taken over the
conformal class of $g$. Then, proceeding in a way analogous to
\cite{ammann.habil,ammann_p04},
one shows that 
$$0 < \la^+ =   \inf \frac{{\left( \int_M |L_g u |^{\frac{2n}{n+2}}\, dv_g
    \right)}^{\frac{n+2}{n}}}{\int_M u L_gu\,dv_g}$$
where the infimum is taken over the smooth functions $u$ such that 
$$\int_M u L_gu\,dv_g>0.$$
Then, one shows using test functions that $\la_+ \leq \mu_1(\mS^n)$. If the
inequality is strict, then we can find a minimizer for the functional above
which is a solution of the Yamabe equation. If the Yamabe invariant is
positive, this solution is a Yamabe metric and hence is positive. However,
if the Yamabe invariant is non-positive, this solution has an alternating
sign.

\section*{A.\ \ Appendix: Proof of Lemma~\ref{regu}}
Let $(M,g)$ be a compact Riemannian manifold of dimension $n \geq 3$ and
let $v \in H_1^2(M)$, $v \not\equiv 0$ and  $u \in L^{N}_+(M)$ be two 
functions which
satisfy in the sense of distributions 
$$L_gv = u^{N-2} v. \eqno{(Eq)}$$
We define $v_+ = \sup(v,0)$. We let $q \in ]1, \frac{n}{n-2}]$ be a fixed
number and $l>0$ be a large real number which will tend to $+\infty$. We
let $\beta= 2q -1 $. We then define the following functions for $x \in
\mR $:

\[ G_l(x) = \left| \begin{array}{ccc}
0 & \hbox{ if } & x <0 \\
x^{\beta}  & \hbox{ if } & x \in
[0,l[\\
l^{q-1}(ql^{q-1} x - (q-1) l^q) 
 & \hbox{ if } & x \geq l
\end{array} \right. \]
and 

\[ F_l(x) = \left| \begin{array}{ccc}
0 & \hbox{ if } & x <0 \\
x^q  & \hbox{ if } & x \in
[0,l[\\
ql^{q-1} x - (q-1) l^q 
 & \hbox{ if } & x \geq l
\end{array} \right. \]

\noindent It is easy to check that for all $x \in \mR$,

\begin{eqnarray} \label{i1}
(F_l'(x))^2 \leq q G_l'(x),
\end{eqnarray}

\begin{eqnarray} \label{i2}
(F_l(x))^2 \geq x G_l(x)
\end{eqnarray}
and 

\begin{eqnarray} \label{i3}
x G'(x) \leq \beta  G_l(x).
\end{eqnarray}

\noindent Since $F_l$ and $G_l$ are uniformly lipschitz continuous
functions, $F_l(v_+)$ and $G_l(v_+)$ belong to $H_1^2(M)$. 
Now, let $x_0 \in M$ be any point of $M$. We denote by $\eta$ a
$C^2$ non-negative function supported in $B_{x_0}(2 \delta)$ ($\delta>0$ being
a small number to be fixed) such that $0 \leq \eta \leq 1$ and such that 
$\eta(B_{x_0}(\delta)) = \{ 1 \}$.
 Multiply equation $(Eq)$ by $\eta^2 G_l(v_+)$
and  integrate over
$M$. Since the supports of $v_+$ and $G_l(v_+)$ coincide, we get:

\begin{eqnarray} \label{equ1} 
c_n \int_M (\nabla v_+, \nabla \eta^2 G_l(v_+) )dv_g + \int_M S_g v_+
\eta^2 G_l(v_+) dv_g = \int_M u^{N-2} v_+  \eta^2 G_l(v_+) dv_g.
\end{eqnarray}

\noindent Let us deal with the first term of the left hand side of
(\ref{equ1}). In the following, $C$ will denote a positive constant
depending possibly on $\eta, q, \beta, \delta$ but not on $l$. We have 
\begin{eqnarray*} 
\int_M (\nabla v_+, \nabla \eta^2 G_l(v_+) )dv_g  & = & 
 \int_M G_l(v_+)  (\nabla v_+,\nabla \eta^2)dv_g +  
  \int_M G_l' (v_+) \eta^2 |\nabla v_+|^2 dv_g \\
& = & \int_M G_l(v_+) v_+ \Delta (\eta^2) - 2 \int_M v_+ G_l'(v_+) \eta 
(\nabla v_+,\nabla \eta) dv_g + \int_M G_l' (v_+) \eta^2 |\nabla v_+|^2 dv_g \\
& \geq &  - C \int_M v_+ G_l(v_+)dv_g  - 2 \int_M v_+^2 G_l'(v_+) |\nabla
\eta|^2 dv_g + \frac{1}{2}   \int_M G_l' (v_+) \eta^2 |\nabla v_+|^2
dv_g.\\
\end{eqnarray*}

\noindent Using (\ref{i1}), (\ref{i2}) and (\ref{i3}),  we get 

\begin{eqnarray}
\int_M (\nabla v_+, \nabla \eta^2 G_l(v_+) )dv_g  & \geq  & 
- C \int_M (F_l(v_+))^2 dv_g + \frac{1}{2q} \int_M (F_l' (v_+))^2 
\eta^2 |\nabla v_+|^2 dv_g  \nonumber \\
& \geq&   -C \int_M (F_l(v_+))^2 dv_g + \frac{1}{2q}\int_M
\eta^2 |\nabla F_l(v_+)|^2 dv_g \nonumber \\
 & \geq & -C \int_M (F_l(v_+))^2 dv_g + \frac{1}{4q} \int_M |\nabla (\eta
 F(v_+)) |^2 dv_g - \frac{1}{2q}  \int_M |\nabla \eta|^2  (F_l(v_+))^2
 dv_g \nonumber \\
& \geq &  -C \int_M (F_l(v_+))^2 dv_g  + \frac{1}{4q} \int_M |\nabla (\eta
 F(v_+)) |^2 dv_g. \label{i4} 
\end{eqnarray}

\noindent Using the Sobolev embedding $H_1^2(M)$ into $L^N(M)$, there
exists a constant $A>0$ depending only on $(M,g)$ such that  
$$ \int_M |\nabla (\eta
 F(v_+)) |^2 dv_g   \geq  A {\left( \int_M (\eta F(v_+))^N dv_g
   \right)}^{\frac{2}{N}} -  \int_M (\eta F(v_+))^2 dv_g.$$
Together with  (\ref{i4}), we obtain 
\begin{eqnarray} \label{i5}
\int_M (\nabla v_+, \nabla \eta^2 G_l(v_+) )dv_g \geq -C \int_M
(F_l(v_+))^2 dv_g  +\frac{A}{4q}  {\left( \int_M (\eta F(v_+))^N dv_g
   \right)}^{\frac{2}{N}} 
\end{eqnarray}

\noindent Independently, we choose $\delta>0$ small enough such that   
$$\int_{B_{x_0}(2 \delta)} u^N dv_g 
\leq {\left( c_n \frac{A}{8q} \right)}^{\frac{n}{2}}.$$
Relation (\ref{i2}) and H\"older inequality  then lead to 
\begin{eqnarray} \label{i6}
\int_M u^{N-2} v_+ \eta^2  G_l(v_+) dv_g \leq 
\int_M u^{N-2} \eta^2 (F_l(v_+))^2 dv_g \leq 
c_n \frac{A}{8q}  {\left( \int_M (\eta F(v_+))^N dv_g
   \right)}^{\frac{2}{N}}.
\end{eqnarray}

\noindent Since, by (\ref{i2}), 
$$\int_M S_g v_+
\eta^2 G_l(v_+) dv_g \geq -C \int_M (F_l(v_+)^2) dv_g,$$
we get from (\ref{equ1}), (\ref{i5}) and (\ref{i6})  that 
$$c_n \frac{A}{8q}  {\left( \int_M (\eta F(v_+))^N dv_g
   \right)}^{\frac{2}{N}} \leq C \int_M
(F_l(v_+))^2 dv_g. $$
Now, by Sobolew embedding, $v_+ \in L^N(M)$. 
Since $2q \leq N$ and since $C$ does not depend on $l$, 
the right hand side of this inequality is
bounded when $l$ tends to $+\infty$. We obtain that  
$$\limsup_{l \to +\infty} \int_M (\eta F(v_+))^N dv_g < +\infty.$$
This proves that $v_+ \in L^{qN}(B_{x_0}(\delta))$. Since $x_0$ is
arbitrary, we get that  $v_+ \in L^{qN}(M)$. Doing the same with
$\sup(-v,0)$ instead of $v_+$, we get that $v \in L^{qN}(M)$. This proves
Lemma \ref{regu}.

\vspace{1cm}               
Authors' address:               
\nopagebreak   
\vspace{5mm}\\   
\parskip0ex     
\vtop{   
\hsize=8cm\noindent   
\obeylines               
Bernd Ammann and Emmanuel Humbert               
Institut \'Elie Cartan BP 239              
Universit\'e de Nancy 1               
54506 Vandoeuvre-l\`es -Nancy Cedex               
France                           
\vspace{0.5cm}               
               
E-Mail:               
{\tt bernd.ammann at gmx.net and humbert at iecn.u-nancy.fr}               
}   
                    
\end{document}